\documentclass[epsf]{article}
\usepackage{amsmath,amsthm,amsfonts,latexsym,amscd}
\usepackage{graphicx}
\title{Minimal two-spheres of low index in manifolds with positive complex sectional curvature}
\author{John Douglas Moore and Robert Ream\\Department of 
Mathematics\\University of California\\
Santa Barbara, CA, USA 93106}

\begin{document}

\maketitle

\begin{abstract}
Suppose that $S^n$ is given a generic Riemannian metric with sectional curvatures which satisfy a suitable pinching condition formulated in terms of complex sectional curvatures.  This pinching condition is satisfied by manifolds whose real sectional curvatures $K_r(\sigma )$ satisfy
$$1/2 < K_r(\sigma ) \leq 1.$$
Then the number of minimal two spheres of Morse index $\lambda $, for $n-2 \leq \lambda  \leq 2n-5$, is at least $p_{3}(\lambda -n+2)$, where $p_{3}(k)$ is the number of $k$-cells in the Schubert cell decomposition for $G_3({\mathbb R}^{n+1})$. \end{abstract}

\section{Introduction}
\label{S:introduction}

Suppose that $M$ is a compact Riemannian manifold, $\hbox{Map}(S^2 ,M)$ the space of smooth maps from the Riemann sphere $S^2$ to $M$.  An element $f \in \hbox{Map}(S^2 ,M)$ is {\em harmonic\/} if it is critical point for the energy function
$$E : \hbox{Map}(S^2 ,M) \longrightarrow {\mathbb R},$$
which is defined by the Dirichlet integral
$$E(f) = \frac{1}{2} \int _{S^2} |df|^2 dA = \frac{1}{2} \int _{S^2} \left[ \left| \frac{\partial f}{\partial x}\right| ^2 +  \left|\frac{\partial f}{\partial y}\right| ^2 \right] dx dy.$$
In the middle integrand of this formula, $|df|$ and $dA$ are calculated with respect to the standard Riemannian metric of constant curvature and total area one on $S^2$, while in the right-hand integrand $(x,y)$ are isothermal coordinates on $S^2$ with respect to the standard conformal structure on the Riemann sphere $S^2$, the integrand being independent of the choice.  If the Riemannian metric on $M$ is induced by an imbedding of $M \subseteq {\mathbb R}^N$ into Euclidean space (always possible by the Nash imbedding theorem), the Euler-Lagrange equation in isothermal parameters for the energy function is
$$\left(\frac{\partial ^2f}{\partial x^2} + \frac{\partial ^2f}{\partial y^2}\right)^\top = 0,$$
where $(\cdot )^\top$ represents orthogonal projection from ${\mathbb R}^N$ into the tangent space.  It is well-known that harmonic maps from the Riemann sphere are also critical points for the area function, and hence also represent parametrized {\em minimal\/} two-spheres, possibly with branch points.

Parametric minimal two-spheres are part of a more general theory which treats parametrized minimal surfaces of arbitrary genus, and is described in the recent book \cite{Mo2}.  The case of minimal two-spheres is simpler in some ways than the higher genus case because there is only one conformal structure on $S^2$, but more complicated in others because it has a positive-dimensional symmetry group of linear fractional transformations, which is noncompact.   

If we regard $S^2$ as the one-point compactification ${\mathbb C} \cup \{ \infty \}$ of the complex plane ${\mathbb C}$ with the standard coordinate $z = x + iy$, and if $T : S^2 \rightarrow S^2$ is the linear fractional transformation (a conformal diffeomorphism) defined by
$$T (z) = \frac{az+b}{cz+d}, \quad \hbox{for} \quad \begin{pmatrix} a & b \cr c & d \end{pmatrix} \in SL(2,{\mathbb C}),$$
then $E(f \circ T) = E(f)$.  Note that $-I \in SL(2,{\mathbb C})$ acts as the identity, so the group of linear fractional transformations is actually
$${\mathcal G} = PSL(2,{\mathbb C}) = \frac{SL(2,{\mathbb C})}{\pm I},$$
and if we want to allow orientation-reversing automorphisms of $S^2$ we must replace $PSL(2,{\mathbb C})$ by the larger group
$$\widehat {\mathcal G} = PSL(2,{\mathbb C}) \cup R \cdot PSL(2,{\mathbb C}),$$
where $R$ is some anticonformal reflection through a great circle.  A nonconstant minimal two-sphere is {\em prime\/} if it is not a branched cover of a minimal two-sphere of smaller area, and two prime minimal surfaces are {\em geometrically distinct\/} if they do not lie on the same $\widehat {\mathcal G}$-orbit.

Our goal here is to develop the Morse theory techniques needed to estimate the number of geometrically distinct prime minimal two-spheres of low Morse index, using a perturbation of the energy proposed by Sacks and Uhlenbeck \cite{SU1}.  This perturbation is the {\em $\alpha $-energy\/}, for $\alpha > 1$, and it is defined by
\begin{equation}E_{\alpha } : \hbox{Map}(S^2 ,M) \longrightarrow {\mathbb R}, \quad E_{\alpha } (f) = \frac{1}{2} \int _{S^2} [(1 + |df|^2)^\alpha - 1]dA, \label{E:alphaenergy} \end{equation}
where $|df|$ and $dA$ are again defined in terms of the canonical Riemannian metric of constant curvature and total area one on $S^2$, there being exactly one such metric up to diffeomorphisms isotopic to the identity.  Note that $E_\alpha (f)$ is strictly increasing as a function of $\alpha $ and $E(f) \leq E_\alpha (f)$.

As explained by Sacks and Uhlenbeck, the $\alpha $-energy satisfies Condition C of Palais and Smale on the completion $L^{2\alpha }_1(S^2,M)$ of $\hbox{Map}(S^2 ,M)$ with respect to the $L^{2\alpha }_1$ Sobolev norm, and $L^{2\alpha }_1(S^2,M) \subseteq C^0(S^2,M)$ when $\alpha > 1$.  These facts make it possible to develop a Liusternik-Schnirelman theory for $E_\alpha $, and a Morse theory for further perturbations of $E_\alpha $ which have Morse nondegenerate critical points.  Note that $E_\alpha (f) \rightarrow E(f)$ as $\alpha \rightarrow 1$, and our goal is to investigate the limit of the critical locus of $E_\alpha $ as $\alpha \rightarrow 1$.

Indeed, we would like to develop a partial Morse theory for the energy $E$ itself, and the first step towards such a Morse theory is the Bumpy Metric Theorem presented in Chapter~5 of \cite{Mo2} (the proof being a revision of an argument from an earlier paper \cite{Mo1.5}). To state this theorem we need some definitions.  By a {\em generic choice of Riemannian metric\/} on $M$ we mean a metric belonging to a countable intersection of open dense subsets of the spaces of $L^2_k$ Riemannian metrics on $M$, as $k$ ranges over the positive integers.  By a {\em nondegenerate critical submanifold\/} for $E : \hbox{Map}(S^2 ,M) \rightarrow {\mathbb R}$ we mean a finite-dimensional submanfold $S \subset \hbox{Map}(S^2 ,M)$ consisting entirely of critical points for $E$ such that the tangent space to $S$ at a given critical point is exactly the null space for the second variation for $E$.  This definition of nondegenerate critical submanifold has been widely used by Bott, and plays a crucial role in his celebrated survey \cite{Bo} of Morse theory.

Once these preparations are made, we can state the Bumpy Metric Theorem (Theorem~5.1.1 of \cite{Mo2}):  For generic choice of Riemannian metric on a manifold $M$ of dimension at least three, all prime compact oriented parametrized minimal two-spheres $f: S^2 \rightarrow M$ are free of branch points and lie on nondegenerate critical submanifolds, each being an orbit for the action of $\widehat {\mathcal G}$.  

To treat Morse theory by perturbation from $E$ to $E_\alpha $, we need to address two key problems:

\vskip .1in
\noindent
{\bf Bubbling.}  As observed by Sacks and Uhlenbeck, a subsequence of $\alpha $-energy critical points will not necessarily converge to a single minimal surface as $\alpha \rightarrow 1$, but instead may converge to a bubble tree, consisting of several minimal two-spheres connected by necks.

\vskip .1in
\noindent
{\bf Covers.}  Branched covers of degree $d$ of a given prime minimal two-sphere form a critical submanifold of real dimension $4d + 2$.  However, in general the space of branched covers of a given degree is not a nondegenerate critical submanifold, so it is not straightforward to calculate its contribution to the topology of $\hbox{Map}(S^2,M)$.  In addition, we need to rule out unbranched covers of nonorientable minimal projective planes.

\vskip .1in
\noindent
The first problem, bubbling, prevents full Morse inequalities from holding for $E$, while the second is much more complicated than multiple covers of closed geodesics, but when the two difficulties can be controlled, one can expect theorems which should be quite similar to those which have been proven for closed geodesics in Riemannian manifolds.

For example, within the theory of closed geodesics, a theorem of Ballmann, Thorbergsson and Ziller \cite{BTZ} states that a generic simply connected compact Riemannian manifold whose real sectional curvatures satisfy the inequality $1/4 < K(\sigma ) \leq 1$ must contain at least $p_2(\lambda - n + 1)$ geometrically distinct smooth closed geodesics of Morse index $\lambda $, where $p_2(k)$ is the number of $k$-cells in the Schubert cell decomposition of $G_2({\mathbb R}^{n+1})$.  Their result can be obtained by establishing equivariant Morse inequalities for the action function
$$J : \hbox{Map}(S^1,M) \rightarrow {\mathbb R},$$
where $\hbox{Map}(S^1,M)$ is the space of loops in $M$.  (The reader can refer to Milnor and Stasheff \cite{MS} for a discussion of the Schubert CW decomposition of $G_m({\mathbb R}^N)$.)

We will prove an analog of this theorem (assuming a stronger pinching constant) for critical points of the energy function $E$ when the Morse index is sufficiently small:

\vskip .1 in
\noindent
{\bf Theorem~1.1.} {\sl Suppose that the compact smooth simply connected manifold $M$ has dimension $n \geq 4$, and has a generic Riemannian metric whose real sectional curvatures $K_r(\sigma )$ satisfy the pinching condition
\begin{equation} \frac{1}{2} < K_r(\sigma ) \leq 1. \label{E:Kr}\end{equation}
Then the number of geometrically distinct prime minimal two spheres of Morse index $\lambda $ is at least $p_{3}(\lambda - n + 2)$ for $n-2 \leq \lambda \leq 2n-5$, where $p_{3}(k)$ is the number of $k$-cells in the Schubert cell decomposition for $G_3({\mathbb R}^{n+1})$.}
  
\vskip .1 in
\noindent
Note that it follows from the differentiable sphere theorem of Brendle and Schoen \cite{BS1} that the pinching condition (\ref{E:Kr}) implies that $M$ is diffeomorphic to an $n$-sphere.

In fact, we will weaken the pinching condition (\ref{E:Kr}) somewhat by making use of complex sectional curvatures.  We recall that in normal coordinates $(x^1 , \ldots , x^n)$ centered at a point $p$ on $M$, the Riemannian metric can be expressed by a Taylor series
$$g_{ij} = \delta _{ij} - \frac{1 }{3} \sum _{k,l} R_{ikjl}(p) x^kx^l + \hbox{ higher order terms,}$$
and the curvature operator can be defined as the linear map ${\mathcal R} : \Lambda ^2T_pM \longrightarrow \Lambda ^2T_pM$ such that
$${\mathcal R}\left(\left. \frac{\partial }{\partial x^i}\right|_p \wedge \left. \frac{\partial}{\partial x^j}\right|_p \right) = \sum _{k,l} R_{ijkl}(p) \left. \frac{\partial}{\partial x^k}\right|_p \wedge \left. \frac{\partial }{\partial x^l}\right|_p.$$
If $z$ and $w$ are linearly independent elements of $T_pM \otimes {\mathbb C}$, the {\em complex sectional curvature\/} of the two-plane $\sigma $ spanned by $z$ and $w$ is
\begin{equation} K(\sigma ) = \frac{\langle {\mathcal R}(z \wedge w), \bar z \wedge \bar w \rangle}{\langle z \wedge w , \bar z \wedge \bar w \rangle }, \label{E:complexsectional}\end{equation}
where the curvature operator ${\mathcal R}$ and the inner product have been extended complex linearly and the bar denotes complex conjugation.  The complex sectional curvature arises quite naturally in the formula for second variation of energy for a minimal two-sphere $f$, as presented in \cite{MM}.  This second variation is the symmetric bilinear form
$$d^2E(f) : T_f\hbox{Map}(S^2,M) \otimes {\mathbb C} \times T_f\hbox{Map}(S^2,M) \otimes {\mathbb C} \longrightarrow {\mathbb C}$$
given by
\begin{equation} d^2E(f)(V,\bar V) =  4\int _{S^2} \left[ \left| \frac{DV}{\partial \bar z} \right|^2 - \left\langle {\mathcal R}\left(V\wedge \frac{\partial f} {\partial z}\right), \bar V \wedge \frac{\partial f} {\partial \bar z} \right\rangle \right] dx dy,\label{E:secondvariationa}\end{equation}
in which
$$\frac{\partial f} {\partial z} = \frac{1}{2} \left( \frac{\partial f} {\partial x} - i \frac{\partial f} {\partial y}\right) \quad \hbox{and} \quad \frac{D}{\partial \bar z} = \frac{1}{2} \left( \frac{D}{\partial x} + i \frac{D}{\partial y}\right)$$
is the Cauchy-Riemann operator defined by the Levi-Civita connection on ${\bf E} = f^*TM \otimes {\mathbb C}$.  If $\sigma $ is the two-plane spanned by $\partial f/\partial z$ and $V$, then (\ref{E:secondvariationa}) can also be written as
\begin{equation}d^2E(f)(V,\bar V) = 2 \int _{S^2} \left[ \left\| \bar \partial V \right\|^2 - K(\sigma ) e(f) \left\| V \right\|^2 \right] dA,\label{E:secondvariation0}\end{equation}
where $K(\sigma )$ is the complex sectional curvature of $\sigma $, $e(f)$ is the energy density of $f$ and $dA$ is the element of area on $S^2$ when the $S^2$ is given the canonical Riemannian metric of constant curvature and total area one.

A theorem of Koszul and Malgrange \cite{KoMa}  gives ${\bf E}$ a holomorphic structure over $S^2$ which makes $V$ a holomorphic section of ${\bf E}$ exactly when $\bar \partial V = 0$ and the Riemann-Roch Theorem from Riemann surface theory then allows us to estimate the dimension of the space of solutions to the equation $\bar \partial V = 0$.  If the complex sectional curvature is positive, this enables us to give a lower bound on the {\em Morse index\/} of $f$, which by definition is the maximum dimension of a complex linear space of sections of ${\bf E} = f^*TM \otimes {\mathbb C}$ on which $d^2E(f)$ is negative definite.

\vskip .1 in
\noindent
{\bf Definition.} A two-plane $\sigma \subseteq T_pM \otimes {\mathbb C}$ is {\em isotropic\/} if it is spanned by linearly independent vectors $v$ and $w$ such that
$$\langle v,v \rangle = \langle v, w \rangle = \langle w, w \rangle = 0.$$
It is {\em half-isotropic\/} if it is spanned by linearly independent vectors $v$ and $w$ such that
$$\langle v,v \rangle = \langle v, w \rangle = 0.$$
Since $v=x+iy$ is isotropic, $x$ and $y$ are orthogonal and thus span a real two plane $\widehat\sigma$ which we call a real two-plane {\em associated\/} to the complex two-plane $\sigma $. 

\vskip .1in
\noindent
Sectional curvature of half-isotropic two-planes give estimates on the Morse index of minimal two-spheres by formula (\ref{E:secondvariation0}), just as the real sectional curvatures give estimates on the Morse index of geodesics.  Having positive half-isotropic curvatures is stronger than having positive isotropic curvatures, and the first assumption gives stronger index estimates.  For example, when $M$ has dimension three there are no isotropic two-planes, while positive half-isotropic curvature reduces to positive Ricci curvature.  One can show that a Riemannian manifold $M$ has nonnegative half-isotropic curvature if and only if $M \times {\mathbb R}$ has nonnegative isotropic curvature.

Positive isotropic curvatures is sufficient to prove the homeomorpism sphere theorem of \cite{MM} in spite of the weaker index estimates for minimal two-spheres.  Brendle has shown that a compact simply connected Riemannian manifold satisfying the stronger condition of having positive half-isotropic curvatures must be diffeomorphic to a sphere; see Theorem~5.10 in the survey \cite{BS2}, which includes references to the proof. 

\vskip .1 in
\noindent
{\bf Lemma~1.2.} {\sl Pinching of the real sectional curvatures $K_r(\sigma_0 )$ and the half-isotropic curvatures $K_i(\sigma_1 )$ are related by
\begin{equation} \delta < K_r(\sigma_0 ) \leq 1 \quad \Rightarrow \quad \frac{4}{3}\delta - \frac{1}{3} < K_i(\sigma_1 ) \leq \frac{4}{3} - \frac{\delta }{3}.\label{E:deltaprime} \end{equation}}

\noindent
Indeed if $(e_1,e_2,e_3,e_4)$ are unit-length vectors, perpendicular to each other, at a point $p$, and $\sigma _1$ is the half-isotropic two-plane spanned by $e_1 + i e_2$ and $ae_3 + ibe_4$, where $a,b \in {\mathbb R}^+$, then
\begin{multline*} \left\langle {\mathcal R}((e_1 + i e_2)\wedge (ae_3 + ibe_4)), (e_1 - ie_2) \wedge (ae_3 - ibe_4) \right\rangle \\ = a^2 (R_{1313} + R_{2323}) + b^2(R_{1414} + R_{2424}) + 2abR_{1243}, \end{multline*}
and according to Berger's formula (7) in \cite{Be}, $|R_{1243}| \leq (2/3)(1-\delta)$, so
\begin{multline*} \delta < K_r(\sigma_0 ) \leq 1 \quad \Rightarrow \\ 2 \delta (a^2 + b^2) - \frac{4}{3}ab(1-\delta) < 2(a^2 + b^2)K_i(\sigma _1) \leq 2(a^2 + b^2) + \frac{4}{3}ab(1-\delta),\end{multline*}
which simplifies to the conclusion of the lemma:
$$ \delta - \frac{1}{3}(1-\delta ) < K_i(\sigma_1 ) \leq 1 + \frac{1}{3}(1-\delta).$$

\vskip .1 in
\noindent
The Lemma shows that if the real sectional curvatures are $1/4$-pinched, the half-isotropic curvatures are positive, and that the hypothesis (\ref{E:Kr}) implies that
\begin{equation} \frac{1}{3} < K_i \quad \hbox{and} \quad K_r \leq 1.\label{E:1/3}\end{equation}
It is actually this weaker hypothesis that is needed in the proof of the estimates for the number of prime minimal two-spheres.

Indeed, the Morse index of minimal two-spheres depends only on {\em pointwise\/} pinching, and the curvature hypothesis in Theorem~1.1 can be weakened: one only needs a relation between the complex sectional curvature of a half-isotropic two-plane $\sigma$ and the real sectional curvature of an associated real two-plane $\widehat\sigma$.

\vskip .1 in
\noindent
{\bf Main Theorem~1.3.} {\sl Suppose that the smooth manifold $M$ is diffeomorphic to a sphere of dimension $n \geq 4$ and has a generic Riemannian metric which satisfies the condition
\begin{equation}  K_i(\sigma) > \frac{1}{3} K_r(\widehat\sigma) > 0, \label{E:1/3bis}\end{equation}
where $\widehat\sigma$ is a real two-plane associated to $\sigma $.  Then the number of distinct prime minimal two spheres of Morse index $\lambda $, none of which double cover minimal projective planes, is at least $p_{3}(\lambda - n + 2)$ for $n-2 \leq \lambda \leq 2n-5$, where $p_{3}(k)$ is the number of $k$-cells in the Schubert cell decomposition for $G_3({\mathbb R}^{n+1})$.}

\vskip .1 in
\noindent
The previous remarks show that Main Theorem~1.3 implies Theorem~1.1.

An earlier article \cite{Mo1} proved the existence of ${n+1 \choose 3}$ minimal two-spheres (of approriate Morse index) for a generic metric on an $n$-manifold $M$ diffeomorphic to a sphere with a stronger condition on real sectional curvature.\footnote{The Corollary of \S 1 of \cite{Mo1} also assumes that $M$ has no minimal projective planes, but it turns out that such projective planes do not interfere with the existence of minimal two-spheres, as we explain in Remark 10.1 at the end of this paper.}  This article asks what happens as the pinching condition is relaxed somewhat, and replaced by a pointwise pinching condition.  We expect that bubbling and branched covers begin to interfere with the achievable partial Morse inequalities as the pinching condition becomes weaker.

To prove Theorem~1.3, we proceed by analogy with the closed geodesic problem to develop an analog of the Morse-Schoenberg comparison theorem (see \S 6.2 of \cite{GKM}), which shows that Morse index grows with energy.  Although the estimates obtained are not as sharp as in the geodesic case, they are sufficient to show that any minimal two-sphere $f : S^2 \rightarrow M$ in a manifold satisfying estimate (\ref{E:1/3bis}) must have Morse index at least $n - 2$, where $n = \dim M$.  This enables us to establish the following step which controls bubbling:

\vskip .1 in
\noindent
{\bf Nonbubbling Theorem~1.4.} {\sl Suppose that the $n$-dimensional Riemannian manifold $M$ has positive half-isotropic sectional curvatures, where $n \geq 4$.  If a sequence of $\alpha $-energy critical points approaches a bubble tree containing $k$ nonconstant minimal two-spheres in the limit, then the Morse index of these critical points must eventually be at least $k(n-2)$.}

\vskip .1 in
\noindent
The index estimates are also used to control branched covers.  We show that under the curvature hypothesis (\ref{E:1/3bis}) all branched covers of order $d \geq 3$ must have Morse index at least $2n-4$, and then use equivariant Morse theory with ${\mathbb Z}_2$ coefficients to rule out branched covers of order two, as well as double covers of projective planes.

To put the ingredients together, we use the constant curvature metric on $S^n$ to determine the low-dimensional equivariant topology of $\hbox{Map}(S^2,M)$, and then use Morse inequalities to derive the estimates on number of minimal two-spheres for generic pinched metrics satisfying our pinching conditions.

Here is a more detailed outline of the proof of the Main Theorem~1.3.  In \S~\ref{S:perturbations}, we describe a center of mass technique which allow us to reduce the noncompact symmetry group $\widehat {\mathcal G}$ for $E$ to the compact symmetry group $\widehat G = O(3)$ needed for $E_\alpha $, and in \S~\ref{S:isotropy} we show that critical points for $E_\alpha $ always have finite isotropy.  In \S~\ref{S:MorseBanach} we explain how to extend Uhlenbeck's Morse theory on Banach manifolds to the $G$-equivariant case, and go on to give Morse index estimates in \S~\ref{S:indexestimates}.  We then prove the Nonbubbling Theorem in \S~\ref{S:nonbubbling}.  We study the space of $d$-fold covers of a prime imbedding in \S\ref{S:branchedcovers}, providing an index estimate which rules out the case $d \geq 3$ under our curvature hypothesis, and describing an explicit model of Milnor \cite{Mil3} for treating the case $d = 2$.  We then review equivariant Morse theory and in \S\ref{S:equivariant} describe its calculation for the space of maps $f :S^2 \rightarrow S^n$ into the round sphere of energy strictly less than $8 \pi$.  The Morse-Witten chain complex approach to Morse theory is described in \S\ref{S:morsewitten}, and we show that $G$-equivariant Morse functions are dense by an argument of Wassermann \cite{Wa} extended to infinite-dimensions.  We then have the ingredients necessary to reduce the proof of Theorem~1.3 to a calculation in equivariant cohomology, using techniques of Bott \cite{Bo} and Hingston \cite{Hing}.  The fact that we need only consider two-fold covers allows us to use ${\mathbb Z}_2$-coefficients which simplify the equivariant cohomology calculations.

The authors thank the referee for many helpful suggestions that improved the article considerably.

\section{Perturbations of minimal surfaces}
\label{S:perturbations}

Our strategy is to show that a harmonic two-sphere in $M$ can be perturbed to a family of $\alpha $-harmonic two-spheres, depending on the parameter $\alpha $, for $\alpha \in (1, \alpha _0)$ for some $\alpha _0 > 1$, but we must deal with the technical problem that the group of symmetries changes when the usual energy $E$ is perturbed to the $\alpha $-energy $E_\alpha $ for $\alpha > 1$.  While the usual energy $E$ is invariant under the group
$$\widehat {\mathcal G} = PSL(2,{\mathbb C}) \cup R \cdot PSL(2,{\mathbb C}) , \quad \hbox{$R$ a reflection,}$$
of conformal and anticonformal transformations of the Riemann sphere, when $\alpha > 1$, the $\alpha $-energy
$$E_\alpha : C^2(S^2,M) \rightarrow {\mathbb R} \quad \hbox{defined by} \quad E_\alpha (f) = {1\over 2} \int _{S^2} (1 + |df|^2)^\alpha dA - \frac{1}{2}$$
is invariant only under the smaller group $O(3)$ of isometries.  This smaller group of symmetries is reflected in an additional condition on critical points of $E_\alpha $: they must be parametrized so that a certain \lq\lq center of mass" is zero:

\vskip .1 in
\noindent
{\bf Lemma 2.1.} {\sl Let ${\bf X} : S^2 \rightarrow {\mathbb R}^3$ denote the standard inclusion, and let ${\bf 0}$ denote the origin in ${\mathbb R}^3$.  For each $\alpha \in (1, \infty)$, there is a smooth function $\psi _\alpha : [0,\infty) \rightarrow {\mathbb R}$ such that $\psi _\alpha(0) = 0$, $\psi _\alpha(t), \psi _\alpha'(t) > 0$ for $t > 0$, and if $f$ is a critical point for $E_\alpha $, then
\begin{equation} \int _{S^2} {\bf X} \psi _\alpha (|df|^2) dA = {\bf 0}. \label{E:centerof mass}\end{equation} }

\noindent
In the argument for this lemma, we regard $S^2$ as the one-point compactification ${\mathbb C} \cup \{ \infty \}$ of the complex plane ${\mathbb C}$ with the standard coordinate
\begin{equation} z = x + iy = re^{i\theta } = e^{u + i \theta }, \label{E:coordinatesonS^2}\end{equation}
and give $S^2$ the constant curvature metric expressed in polar coordinates by
\begin{equation} g = ds^2 = \frac{4}{(1+r^2)^2} (dr^2 + r^2 d \theta ^2) = \frac{1}{\cosh ^2u}(du^2 + d \theta ^2), \label{E:constatcurvmetric}\end{equation}
where $r = e^{u}$.  (This metric is related by constant rescaling to the constant curvature metric of total area one.)  For $t \in {\mathbb R}$, we define a family of linear fractional transformations,
$$\phi_t : S^2 \rightarrow S^2 \quad \hbox{by} \quad u \circ \phi_t = u + t, \quad \theta \circ \phi_t = \theta ,$$
so that $t \mapsto \phi_t$ is a one-parameter subgroup of conformal transformations.  The energy density is then given by the formula
$$e(f\circ \phi _t) = \frac{1}{2} \left( \left| \frac{\partial f}{\partial u} \right| ^2 + \left| \frac{\partial f}{\partial \theta } \right| ^2\right) ^2 \cosh^2(u+t),$$
and it is straightforward to calculate its derivative at $t = 0$:
\begin{multline} \left. \frac{d}{dt}e(f \circ \phi _t) \right|_{t=0} = \frac{1}{2} \left. \frac{d}{dt}|d(f \circ \phi _t)|^2 \right|_{t=0} \\ = \left( \left| \frac{\partial f}{\partial u} \right| ^2 + \left| \frac{\partial f}{\partial \theta } \right| ^2 \right) \sinh u \cosh u = |df|^2 \tanh u. \label{E:derivofefphi}\end{multline}
In terms of the standard Euclidean coordinates $(x,y,z)$ on ${\mathbb R}^3$, $S^2$ is represented by the equation $x^2 + y^2 + z^2 = 1$, and a straightforward computation using stereographic projection from the north pole to the $(x,y)$-plane shows that
$$z = \frac{r^2-1}{r^2+1} = \frac{\sinh u}{\cosh u} = \tanh u.$$

We can use (\ref{E:derivofefphi}) to calculate the derivative of $E_\alpha $, and obtain
\begin{multline*} \left. \frac{d}{dt} E_\alpha (f \circ \phi _t) \right|_{t=0} = \frac{1}{2} \left. \frac{d}{dt} \int_{S^2} (1 + |d(f \circ \phi _t)|^2)^\alpha \hbox{sech}^2(u+t) du d\theta \right|_{t=0} \\ = \alpha \int _{S^2} (1 + |df|^2)^{\alpha - 1}\tanh u |df|^2 dA - \int _{S^2} (1 + |df|^2)^\alpha \tanh u \ dA. \end{multline*}
Thus, if $f$ is a critical point for $E_\alpha $,
\begin{multline*} 0 = \int_{S^2} [\alpha  (1 + |df|^2)^{\alpha - 1}  |df|^2 - (1 + |df|^2)^{\alpha} ] zdA, \\ =  \int_{S^2} [ \alpha (1 + |df|^2)^{\alpha - 1}  |df|^2 - (1 + |df|^2)^{\alpha} + 1] zdA,\end{multline*}
where we have used the fact that the average value of $z$ on $S^2$ is zero.  We obtain (\ref{E:centerof mass}) by setting
\begin{equation} \psi _\alpha (t) = \frac{\alpha (1+t)^{\alpha -1} t - (1 + t)^\alpha +1}{\alpha -1} = \int _0^t \alpha \frac{\tau}{(1+ \tau)^{2-\alpha}} d \tau. \label{E:centermassalpha}\end{equation}
Conversely, if the center of mass condition (\ref{E:centerof mass}) holds, $f$ is critical for $E_\alpha $ under conformal dilations.  Moreover, $\psi _\alpha (t)$ has a smooth limit as $\alpha \rightarrow 1$, namely
\begin{equation} \psi _1(t) = \int _0^t \frac{\tau }{(1+ \tau)} d \tau = t - \log(1 + t). \label{E:centermass1}\end{equation}

\vskip .1 in
\noindent
{\bf Definition.} We will say that a smooth map $f : S^2 \rightarrow M$ has {\em $\alpha $-center of mass zero\/} if it satisfies the condition
\begin{equation} \int _{S^2} {\bf X} \psi _\alpha (|df|^2) dA = {\bf 0}, \label{E:alphacenter}\end{equation}
where ${\bf X} : S^2 \rightarrow {\mathbb R}^3$ is the standard inclusion as a round sphere centered at the origin. and we let
$$\hbox{Map}_0(S^2,M) = \left\{ f \in \hbox{Map}(S^2,M) : \int _{S^2} {\bf X} \psi _1(|df|^2) dA = {\bf 0} \right\}.$$

\noindent
Note that if $T \in \widehat {\mathcal G}$,
$$f, \ f \circ T \in \hbox{Map}_0(S^2,M) \quad \hbox{and $f$ nonconstant} \quad \Rightarrow \quad T \in O(3),$$
so setting the center of mass zero does indeed cut the symmetry group down to $O(3)$.

\vskip .1 in
\noindent
{\bf Theorem 2.2.} {\sl Suppose that $M$ is a compact connected Riemannian manifold and that $f : S^2 \rightarrow M$ is a prime oriented parametrized minimal immersion lying on a nondegenerate critical submanifold $N$ for $E$ of dimension six, an orbit for the $PSL(2,{\mathbb C})$ action described above.  Let $N_0 = N \cap \hbox{Map}_0(S^2,M)$.  Then for some $\alpha _0 \in (1,\infty)$, there is a smooth map
$$\eta : N_0 \times [1,\alpha _0) \longrightarrow \left\{ (f, \alpha ) \in L^{2}_{k}(S^2 ,M) \times [1,\alpha _0) : \int _{S^2} {\bf X} \psi _\alpha (|df|^2) dA = {\bf 0} \right\}$$
which satisfies the following conditions:
\begin{enumerate}
\item $\eta (f, \alpha )$ is a critical point for $E_\alpha $ for each $f \in N_0$ and $\alpha \in [1,\alpha _0)$,
\item $\eta $ is $O(3)$-equivariant, and
\item there is an $\epsilon _0 > 0$ such that if $f'$ is a critical point for $E_\alpha $ such that $\| f' - f \|_{C^2} < \epsilon _0$ for some $f \in N_0$, $f'$ lies in the image of $\eta $. 
\end{enumerate}
}

\vskip .1 in
\noindent
Note that it follows from Condition 2 that if $f \circ A = f$ for $f \in N_0$, where $A : S^2 \rightarrow S^2$ is the antipodal map, then 
$$\eta (f, \alpha ) \circ A =  \eta (f, \alpha ), \quad \hbox{for $f \in N_0$ and $\alpha \in [1,\alpha _0)$.}$$

%Reference to $\omega $-lemma included:
Once we have dealt with the change in symmetry group, the proof of Theorem~2.2 is a relatively direct consequence of the implicit function theorem.  Note that when $k$ is relatively large, the map
$$E_\alpha : L^{2}_k(S^2,M) \longrightarrow {\mathbb R}$$
is actually $C^\infty $; indeed, we can regard $E_{\alpha}$ as a composition of several maps
$$f \mapsto df \mapsto (1 + |df|^2)^\alpha \mapsto \frac{1}{2} \int _\Sigma (1 + |df|^2)^\alpha dA - \frac{1}{2},$$
each of which is smooth.  The first map is smooth into the Hilbert manifold $L^2_{k-1}(S^2,TM)$, the second is smooth by the $\omega $-Lemma (see Theorem~1.4.7 of \cite{Mo2}), while the third map is smooth since integration is continuous and linear.  Similarly, the two-variable map
$$E^\star : L^{2}_{k}(S^2 ,M) \times [1,\alpha _0) \longrightarrow {\mathbb R} \quad \hbox{defined by} \quad E^\star (f, \alpha ) = E_\alpha (f),$$
is also smooth when $\alpha _0 > 1$, as is its restriction to
$${\mathcal M}_0 = \left\{ (f, \alpha ) \in L^{2}_{k}(S^2 ,M) \times [1,\alpha _0) : \int _{S^2} {\bf X} \psi _\alpha (|df|^2) dA = {\bf 0} \right\}.$$
Of course, when $\alpha = 1$ we obtain the usual energy, $E^\star (f,1) = E(f)$ restricted to the $L^2_k$-completion of $\hbox{Map}_0(S^2,M)$.

Choose $f_0 \in N_0$ and let $S_{f_0}$ be a slice of the $SO(3)$-action through $(f_0,1) \in {\mathcal M}_0$.  We choose a coordinate system
$$\phi : S_{f_0} \longrightarrow V \times [1, \alpha _0) \subseteq {\bf H} \times [1, \alpha _0),$$
where $V$ is a convex open neighborhood of $0$ in a Hilbert space ${\bf H}$ with inner product $\langle \cdot , \cdot \rangle$, such that
$$\phi (f,\alpha ) = (\phi _0(f, \alpha), \alpha ).$$
Using this coordinate system we create a corresponding smooth function
$$F^\star : V \times [1, \alpha _0) \longrightarrow {\mathbb R}, \quad F^\star(\phi _0(f, \alpha), \alpha ) = E^\star(f,\alpha ).$$
We then define a vector field ${\mathcal X} : V \times [1, \alpha _0) \rightarrow {\bf H}$ by
$$D_1F^\star(x,\alpha) = \langle {\mathcal X}(x),y \rangle, \quad \hbox{for $x,y \in {\bf H}$.}$$
Since the null space for $d^2E(f)$ is generated by the symmetry group $G$, the restriction of $d^2E(f)$ to $S_{f_0}$ is nondegenerate, and this implies that $D_1{\mathcal X}(0)$, the partial derivative with $\alpha = 1$ held constant, is invertible.  The implicit function theorem (eg Theorem 5.9 in Chapter I of \cite{L}) implies that there is a smooth map $g : [1, \alpha _0) \rightarrow {\bf H}$ with $g(1) = 0$ such that
$${\mathcal X}(g(\alpha ), \alpha ) = 0.$$
This in turn yields a smooth curve $\tilde g : [1, \alpha _0) \rightarrow {\mathcal M}_0$, $\tilde g(\alpha ) = (\eta (\alpha ), \alpha )$, in which each $\eta (\alpha )$ is a critical point for the restriction of $E_\alpha $ to the space of maps satisfying (\ref{E:alphacenter}).  But our calculation at the beginning of the section shows that these are also critical for conformal dilations, hence critical for all variations within $L^{2}_{k}(S^2 ,M)$.

Thus, after possibly decreasing $\alpha _0 > 1$, we obtain a smooth curve of critical points
$$\alpha \mapsto \eta(\alpha ) = \eta (f_0, \alpha )$$
for $E_\alpha $.  We then extend $\eta $ to $N_0 \times [1,\alpha _0)$ by demanding $SO(3)$-equivariance, and the resulting extension satisfies all the conditions of Theorem 2.2.

\section{Isotropy groups}
\label{S:isotropy}

When developing general theory for the next few sections, we assume that $M$ is a compact smooth manifold with an arbitrary Riemannian metric.  The $\alpha $-energy
$$E_\alpha : \hbox{Map}(S^2 ,M) \longrightarrow {\mathbb R}$$
is invariant under an action of the group $\widehat G = O(3)$ by isometries on the domain.  Because of this action, the nonconstant critical points for $E_\alpha $ lie on $\widehat G$-orbits of critical points, and our goal is to apply $\widehat G$-equivariant Morse theory, in which isotropy groups of the $\widehat G$-action play an important role.

The isotropy group of any point in $\hbox{Map}(S^2 ,M)$ must be a closed subgroup of $\widehat G$, hence a Lie subgroup.  One possibility it that the isotropy group is $\widehat G$ itself, and this occurs exactly for the constant maps, which form a submanifold $M_0$ of zero energy which is diffeomorphic to $M$ itself.  If the isotropy group $H \subseteq \widehat G$ of a nonconstant map is not discrete, the only possibility for its identity component is $S^1$, a group of rotations about some axis.

If the isotropy group contains $S^1$, a smooth map $f : S^2 \rightarrow M$ must in fact degenerate to a curve parametrized by $S^2$, so we examine whether this could happen for a nonconstant critical point for $E_\alpha $.  As in the previous section, we give $S^2$ the constant curvature metric expressed in polar coordinates by
$$g = ds^2 = \frac{4}{(1+r^2)^2} (dr^2 + r^2 d \theta ^2) = \frac{1}{\cosh ^2u}(du^2 + d \theta ^2),$$
where $r = e^{u}$, chosen so that a rotation through an angle $\phi $ is represented by $\theta \mapsto \theta + \phi$.  Then critical points for $E_\alpha $ which are invariant under $S^1$ are simply critical points independent of $\theta$, so we restrict $E_\alpha $ to functions $f_0(u)$ of the variable $u$.  Then $E_\alpha $ simplifies to
$$E_\alpha (f_0) = \pi \int _{-\infty}^\infty \left[1 + (\cosh  u \| f_0'(u) \| )^{2} \right]^\alpha \hbox{sech}^2u \ du - \frac{1}{2},$$
where $f_0 : {\mathbb R} \rightarrow M$ is now a smooth curve, with tangent vector $f_0'(u)$, at $f_0(u)$, which has length $\| f_0'(u) \|$ in terms of the Riemannian metric on $M$.  For $f_0$ to be a critical point, it must satisfy the Euler-Lagrange equation, which is simply
$$\frac{D^g}{du}\left\{ \left[ 1 + (\cosh u \| f_0'(u) \| )^2 \right]^{\alpha - 1} f_0'(u) \right\} = 0,$$
where $D^g$ is the covariant derivative defined by the Riemannian metric on $M$, and this equation has as solutions $f _0 = \gamma \circ h$, where $\gamma : {\mathbb R} \rightarrow M$ is a unit speed geodesic and $h : {\mathbb R} \rightarrow {\mathbb R}$ is a smooth function such that
$$\left[ 1 + (\cosh u \ h'(u))^2 \right]^{\alpha - 1} h'(u) = a,$$
with $a$ being a constant, which is nonzero since $h'$ is nonzero when $f_0$ is nonconstant.  An important feature of this equation is that $h'(u)$ assumes its minimum value when $u = 0$, and hence
$$h'(u) \geq c \quad \hbox{for all $u$,} \quad \hbox{where} \quad g(c) = \left[ 1 + c^2 \right]^{\alpha - 1} c = a.$$
Since $g$ is a strictly increasing function, there is a unique positive solution $c$ to the equation $g(c) = a$, and
\begin{multline*} E_\alpha (f_0) \geq \pi \int _{-\infty}^\infty \left[1 + c^2 \cosh  ^2u  \right]^\alpha  \hbox{sech}^2u \ du - \frac{1}{2} \\ \geq \pi \int _{-\infty}^\infty c^{2\alpha } \cosh ^{2(\alpha - 1)}u \ du - \frac{1}{2} = \infty. \end{multline*}
Thus we see that there are no nonconstant finite energy critical points to $E_\alpha $ which have isotropy group containing $S^1$, or:

\vskip .1 in
\noindent
{\bf Lemma~3.1.} {\sl Critical points for the $\alpha $-energy $E_\alpha : \hbox{Map}(S^2 ,M) \rightarrow {\mathbb R}$ of finite $\alpha $-energy are either constant with isotropy group $\widehat G = O(3)$, or their isotropy groups for the $\widehat G$-action are finite.}

\vskip .1 in
\noindent
Of course, $\hbox{Map}(S^2 ,M)  - M_0$ contains many points at which the isotropy group of the $\widehat G$-action is a circle, but none of these are critical points for $E_\alpha $.  For a given choice of Riemannian metric on $M$, we let $K$ denote the set of critical points for
$$E_\alpha : \hbox{Map}(S^2 ,M)  - M_0 \longrightarrow {\mathbb R}.$$
Then Lemma~3.1 implies that the action of $\widehat G = O(3)$ on an open neighborhood of $K$ is locally free.  Since the local flow for the gradient of $\alpha $-energy preserves isotropy groups, it is also the case that the $\widehat G$-action is locally free on a neighborhood $U$ of all trajectories for this gradient flow which connect points of $K$.

We can also ask what happens when the action of the smaller group $G = SO(3)$ fails to be free at a critical point $f : S^2 \rightarrow M$ for either $E$ or $E_\alpha $.  By Lemma~3.1, the isotropy group must be finite, and the finite subgroups of $SO(3)$ are classified in \S 2.6 of \cite{Wo}.  The only possibilities are cyclic groups, dihedral groups and symmetry groups for the five platonic solids.  Each of these subgroups contains a nontrivial cyclic subgroup which must act on $S^2$ by rotations about some axis, the poles of which are left fixed.  Thus for each of these groups, the poles of rotations are points on $S^2$ with nontrivial isotropy, and in the case of harmonic maps, such points must be branch points.

Under the assumption that all prime minimal two-spheres of Morse index $\leq 2n-5$ are free of branch points, the isotropy of the $SO(3)$ action at these critical points must be trivial.  Critical points for $E_\alpha $ which closely approximate these minimal two-spheres must also have trivial $SO(3)$ isotropy.

\section{Morse theory on Banach manifolds}
\label{S:MorseBanach}

Suppose that ${\mathcal M}$ is a smooth manifold (say an infinite-dimensional Banach manifold) with a continuous right $G$-action by diffeomorphisms, where $G$ is a compact Lie group.  If $F : {\mathcal M} \rightarrow {\mathbb R}$ is a $G$-invariant $C^2$ function, then $F$ descends to a map $F_0 : {\mathcal M}/G \rightarrow {\mathbb R}$, where ${\mathcal M}/G$ is the space of orbits.  When the action of $G$ is free, the quotient is a well-behaved topological manifold.  If the action is only \lq\lq locally free," that is the isotropy groups are finite, we can apply the slice and tube theorems, extended from \S 2.3 and 2.4 of \cite{DK} to Banach manifolds, and study $F$ by means of its restriction to a slice ${\mathcal S}_{f_0}$ through a given point $f_0 \in {\mathcal M}$.

There are three approaches to constructing a Morse theory for such a function $F$:
\begin{enumerate}
\item We can perturb $F$ to a {\em Morse function\/} which has nondegenerate critical points in the usual sense, but this breaks the symmetry under the $G$-action.
\item We can look for a perturbation of $F$ to an {\em equivariant Morse function\/} for which the critical points lie on nondegenerate critical submanifolds (as defined by Bott), each of which is a single $G$-orbit.  If the action is locally free, the restriction of $F$ to any slice has critical points which are nondegenerate in the usual sense. 
\item We can seek a a {\em Morse-Bott function\/} on ${\mathcal M}$ for which all critical points lie on finite-dimensional compact nondegenerate critical submanifolds (which may be larger than a single $G$-orbit).
\end{enumerate}
All three approaches play a role in the proof of our theorems.  The main thread of our argument applies the second case to the $\alpha $-energy $E_\alpha : {\mathcal M} \rightarrow [0,\infty)$, where ${\mathcal M} = L^{2\alpha }_1(S^2 ,M)$, $M$ being a compact Riemannian manifold with generic metric.  On the other hand, the first case is useful for proving ancillary results; it can be applied to a suitable perturbation of the $\alpha $-energy, say $F : {\mathcal M} \rightarrow {\mathbb R}$, where
\begin{equation} F (f) = {1\over 2} \int _{S^2} [(1 + |df|^2)^\alpha - 1]dA + \int _{S^2} f \cdot \psi dA, \label{E:definitionF}\end{equation}
the last perturbation term being defined in terms of an isometric immersion of $M$ into some Euclidean space ${\mathbb R}^N$ and via a suitable small function $\psi : S^2 \rightarrow {\mathbb R}^N$.  Theorem~4.5.3 of \cite{Mo2} states that for generic choices of $\psi $, the function $F$ defined by (\ref{E:definitionF}) is bounded below and has nondegenerate critical points.  The third case applies to the function $E_\alpha $ restricted to the set on which $E_\alpha < 8 \pi$ for the constant curvature one metric on $S^n$, which we compare to generic metrics on $S^n$.   In all cases the function in question is $C^2$, satisfies Condition C and has $C^\infty $ critical points.  

Following \S 4.5.3 of \cite{Mo2}, we now review the Morse theory on Banach manifolds of Uhlenbeck \cite{U}, together with the straightforward extensions to equivariant Morse functions or Morse-Bott functions.  A calculation shows that the second variation of $E_\alpha $ at a critical point $f_0 \in {\mathcal M}$ is the symmetric bilinear form
$$d^2E_\alpha (f_0) : T_{f_0}L^{2\alpha }_1(S^2 ,M) \times T_{f_0}L^{2\alpha }_1(S^2 ,M) \longrightarrow {\mathbb R}$$
defined by
\begin{multline}d^2E_\alpha (f_0)(X,Y) =  \alpha \int _{S^2} (1 + | df _0|^2)^{\alpha -1} [\langle DX, DY \rangle - \langle {\mathcal K}(X),Y \rangle ]dA \\ + 2 \alpha (\alpha - 1) \int _{S^2} (1 + | df _0|^2)^{\alpha - 2} \langle df_0, DX\rangle \langle df_0 , DY \rangle dA,\label{E:secondvariation1}\end{multline}
where $D$ is the covariant differential and the bundle endomorphism ${\mathcal K}$ is determined at $p$ in terms of an orthonormal basis $(e_1,e_2)$ for $T_pM$ by
$$\langle {\mathcal K}(X),Y\rangle(p) = \sum _{i=1}^2\left\langle {\mathcal R}(X(p) \wedge e_i), Y(p) \wedge e_i\right\rangle .$$
This second variation (\ref{E:secondvariation1}) approaches the second variation of the usual energy $E$ in real form
\begin{equation} d^2E(f_0)(X,Y) =  \alpha \int _{S^2} [\langle DX, DY \rangle - \langle {\mathcal K}(X),Y \rangle ]dA, \label{E:realform}\end{equation} 
if $f_0$ approaches a harmonic map as $\alpha \rightarrow 1$.  To see this, one uses the estimate
\begin{equation} \frac{\langle df_0, DX\rangle \langle df_0 , DX \rangle }{1 + | df _0|^2} \leq \langle DX,DX\rangle \label{E:quotientestimate}\end{equation}
to show that the second term in (\ref{E:secondvariation1}) is dominated by
$$\alpha \int _{S^2} (1 + | df _0|^2)^{\alpha -1} [\langle DX, DY \rangle] dA$$
when $\alpha $ is close to one.  Similarly, the second variation of $F$ is given by
\begin{multline}d^2F (f_0)(X,Y) =   \alpha \int _{S^2} (1 + | df_0 |^2)^{\alpha -1} [\langle DX , DY\rangle - \langle {\mathcal K}(X),Y \rangle ]dA \\ + 2 \alpha (\alpha - 1) \int _{S^2} (1 + | df _0|^2)^{\alpha - 2} \langle df_0, DX\rangle \langle df_0 , DY \rangle dA + \int_{S^2} \psi \cdot \alpha (X,Y) dA, \label{E:secondvariation2}\end{multline}
where $\alpha $ in the last term is the second fundamental form of $M$ in ${\mathbb R}^N$,

We define the {\em Morse index\/} of $F$ at an element $f_0 \in L^{2\alpha }_1(S^2,M)$ (whether $f_0$ is critical or not) to be the maximal dimension of a linear subspace of
$$T_{f_0}L^{2\alpha }_1(S^2 ,M)$$
on which $d^2F$ is negative-definite; note that the Morse index is a lower semicontinuous function of $f_0 \in {\mathcal M}$.  It will be useful to compare these formulae to an inner product $\langle\langle \cdot , \cdot \rangle\rangle$ we can define by
\begin{equation} \langle\langle X, Y \rangle\rangle =  2 \alpha \int _{S^2} (1 + | df_0 |^2)^{\alpha -1} [\langle DX , DY\rangle + \langle X,Y \rangle ]dA, \label{E:innerproduct} \end{equation}
which is equivalent to the usual $L^2_1$ inner product, when $f_0$ is a fixed $C^\infty $ map.

It is convenient to extend the notion of \lq\lq gradient-like" vector field ${\mathcal X}$ (as described in Milnor \cite{Mil2}) to the three cases.  In each case, the gradient-like vector field ${\mathcal X}$ is constructed by an explicit formula near the critical locus, which is then pieced together with a pseudogradient as described by Palais \cite{Pal1} away from the critical locus.

\vskip .1in
\noindent
{\bf I. Near the critical locus.}  Consider first the function $F : {\mathcal M} \rightarrow [0,\infty)$ defined by (\ref{E:definitionF}) which is Morse nondegenerate in the usual sense; the critical points of such a function are isolated.  Our first requirement is that the restriction of ${\mathcal X}$ to an open neighborhood $U_{f_0} = U$ of a critical point $f_0 \in {\mathcal M}$ is replaced by a \lq\lq canonical Fredholm map" $A$ which is transported via an exponential chart to $U$,
\begin{equation} A : T_{f_0}{\mathcal M} \longrightarrow T_{f_0}{\mathcal M} \quad \mapsto \quad {\mathcal X}|U = {\mathcal X}_A = \hbox{exp}_{f_0} \circ A \circ (\hbox{exp}_{f_0})^{-1}|U. \label{E:0}\end{equation}
This canonical Fredholm map $A$ is defined in terms of the second variation (\ref{E:secondvariation2}) and the inner product (\ref{E:innerproduct}) by
\begin{equation}d^2F(f_0)(A(X),Y) = \langle\langle X,Y \rangle\rangle. \label{E:canonicalfredholm} \end{equation}
(Thus $A^{-1}$ would define a gradient with respect to the inner product.)  It is straightforward to see that $A$ is a bounded symmetric endomorphism in terms of the $L^2_1$ norm and that it is also bounded in terms of the $L^2_k$ norms for $k \in {\mathbb N}$.  Lemma~5.3 of \cite{U} shows that $A$ restricts to a bounded endomorphism of $L^{2\alpha }_1(S^2 , f_0^*TM)$, while other estimates in \cite{U} show that ${\mathcal X}_A$ satisfies the condition
\begin{equation}dF({\mathcal X}_A)(p) > \varepsilon _0 \| p \|^2 > 0 \quad \hbox{for $p \in U - \{ f_0 \}$,} \label{E:0.5}\end{equation}
for some $\varepsilon _0 > 0$, where $\| \cdot \|$ denotes the $L^{2\alpha }_1$ norm in exponential coordinates.  This is the condition used by Uhlenbeck \cite{U} for the definition of \lq\lq weak gradient," and is sufficient for constructing a handle-body decomposition.

For the other cases, cases in which
$$E_\alpha : {\mathcal M} \rightarrow {\mathbb R}, \quad {\mathcal M} = L^{2\alpha }_1(S^2 ,M)$$
is equivariant Morse or Morse-Bott, the critical locus is the union of compact nondegenerate critical submanifolds.  If $K$ is one of those nondegenerate critical submanifolds and $f_0 \in K$, we can use the inner product (\ref{E:innerproduct}) to divide $T_{f_0}{\mathcal M}$ into an orthogonal direct sum decomposition,
$$T_{f_0}{\mathcal M} = T_{f_0}K \oplus S_{f_0},$$
with $T_{f_0}K$ being not only the tangent space to $K$, but also the null space of the symmetric bilinear form (\ref{E:secondvariation1}).  As $f_0$ varies through $K$, we get a direct sum decomposition
$$T{\mathcal M}|K \oplus S,$$
where $S$ is the normal bundle to $K$ in ${\mathcal M}$.  For any $f_0 \in K$, we can then define $A_{f_0} : S_{f_0} \rightarrow S_{f_0}$ by (\ref{E:canonicalfredholm}), a linear map which is bounded with respect to the $L^{2\alpha }_1$ norm or any of the $L^2_k$ norms restricted to $S_{f_0}$.

The family of maps $f_0 \rightarrow A_{f_0}$ defines a vector field $\tilde {\mathcal X}_A$ on the total space of $S$ which is tangent to the fibers.  Finally, there is a diffeomorphism $\phi _K$ from a neighborhood $V$ of the zero section in this normal bundle to an open neighborhood $U = U_K$ of $K$ such that $\phi _K$ and $(\phi _K)_*$ are the identity along the zero section, and we set ${\mathcal X}_A = (\phi _K)_*(\tilde {\mathcal X}_A)$.  This field satisfies
\begin{equation}dF({\mathcal X}_A)(p) > \varepsilon _0 \| v \|^2 > 0 \quad \hbox{ for } p = \phi _K(v), v \in S_{f_0}, \label{E:0.5a}\end{equation}
for some $\varepsilon _0 > 0$, where $\| \cdot \|$ denotes the $L^{2\alpha }_1$ norm on $S_{f_0}$.

Note that the Morse-Bott case contains the others as special cases.  

\vskip .1in
\noindent
{\bf II. Away from the critical locus.}  The second requirement is that away from the critical locus, the vector field is a {\em pseudogradient\/} as described by Palais \cite{Pal1}; that is, there exist constants $\varepsilon _1 > 0$ and $\varepsilon _2 > 0$ such that 
\begin{equation} \| {\mathcal X}(p) \| < \varepsilon _1  \|dF_p \|, \quad \|dF_p\|^2 < \varepsilon _2 dF_p({\mathcal X}(p)). \label{E:pseudogradient}\end{equation}
When these estimates hold for ${\mathcal X}_0$ and ${\mathcal X}_1$, they hold for any convex combination thereof, and therefore pseudogradients can be pieced together with the $C^2$ partitions of unity available on the Banach manifold $L^{2\alpha }_1(S^2 , f_0^*TM)$.  (That such $C^2$ partitions of unity exist for $L^{2\alpha }_1(S^2 , f_0^*TM)$ is proven near the end of \S 1.11 in the first author's book \cite{Mo2}, this being a special case of the theory of partitions of unity on Banach manifolds described in \S 2C of the survey by Eells \cite{Ee}.)  These estimates imply a third,
\begin{equation} \| {\mathcal X}(p) \|^2 < \varepsilon _1^2 \varepsilon _2 dF_p(X(p)). \label{E:2}\end{equation}
which in turn implies that the flow for $-{\mathcal X}$ will be positively complete by the following argument:

Let $\{ \phi _t : t \in {\mathbb R} \}$ denote the local one-parameter group generated by $-{\mathcal X}$, and let $d :{\mathcal M} \times {\mathcal M} \rightarrow {\mathbb R}$ denote the distance function defined by the Finsler metric which ${\mathcal M}$ inherits from the ambient Finsler manifold $L^{2\alpha }_1(S^2,{\mathbb R}^N)$.  Then it follows from (\ref{E:2}) and the Cauchy-Schwarz inequality that if $t_1 < t_2$,
\begin{multline} d(\phi _{t_1}(p) ,\phi _{t_2}(p))^2 \leq \left[ \int _{t_1}^{t_2} \left\| \frac{d}{dt} \phi _t(p) \right\| dt \right]^2 = \left[ \int _{t_1}^{t_2} \left\| - {\mathcal X}(\phi _t(p)) \right\| dt \right]^2 \\
 \leq (t_2-t_1) \left[ \int _{t_1}^{t_2} \| {\mathcal X}(\phi _t(p)) \|^2 dt \right] \leq \varepsilon _1^2 \varepsilon _2 (t_2-t_1) \int _{t_1}^{t_2} dF({\mathcal X})(\phi _t(p))dt \\
= \varepsilon _1^2 \varepsilon _2 (t_2-t_1) (F(\phi _{t_1}(p)) - F(\phi _{t_2}(p))). \label{E:justofdeformation} \end{multline}
Since $F$ is bounded below, this inequality implies that if $\phi _t(p)$ is defined only for $t \in [0, b)$, for some finite $b$, then any sequence $t_i \rightarrow b$ yields a Cauchy sequence $\{ \phi _{t_i}(p) \}$, and completeness of the Finsler metric shows that this Cauchy sequence must converge to some $q \in {\mathcal M}$.  But there is a family of integral curves passing through points near $q$, showing that $[0, b)$ is not a maximal interval on which $\phi _t(p)$ is defined.

On the other hand, the second inequality of (\ref{E:pseudogradient}) shows that $\|dF \| \rightarrow 0$ along a sequence of points on each orbit, and then Condition C implies that each orbit eventually enters a neighborhood of some critical point.  Since the function $E_\alpha $ or $F$ satisfies Condition C, each component of the critical set $K$ will be compact. 

\vskip .1 in
\noindent
{\bf Definition.}   Suppose that $F: {\mathcal M} \rightarrow {\mathbb R}$ is a $C^2$ Morse-Bott function on the Banach manifold ${\mathcal M}$ which satisfies Condition C and is bounded below.  A {\em gradient-like vector field} for $F$ is a vector field ${\mathcal X}$ on ${\mathcal M}$ which satisfies the following requirements:  If $K_0$ is a component of the critical locus and thus a compact nondegenerate critical submanifold, there exist open neighborhoods
$$V(K_0) = \bigcup \{ V_{f_0} : f_0 \in K_0 \} \quad \hbox{and} \quad U(K_0) = \bigcup \{ U_{f_0} : f_0 \in K_0 \}$$
where $V_{f_0}$ is an open ball neighborhood of $f_0$ in terms of exponential coordinates, with $\overline{V(K_0)} \subseteq U(K_0)$ such that ${\mathcal X}$ satisfies (\ref{E:0}) and (\ref{E:0.5}) on $U_{f_0}$ for each $f_0 \in K_0$.  Moreover, ${\mathcal X}$ is a pseudogradient on
$${\mathcal M} - \bigcup \{ \overline{V(K_0)} : \hbox{ $K_0$ a component of the critical locus } \}.$$

\vskip .1 in
\noindent
A gradient-like vector field is also positively complete by the same argument we gave for pseudogradients, except we replace estimate (\ref{E:2}) by a combination of estimates (\ref{E:0.5}) and (\ref{E:2}).  By Condition C, there are only finitely many critical submanifolds with values of $F$ no larger than a given energy bound $E_0$, and one checks that any orbit for $-{\mathcal X}$ can enter neighborhoods $V(K_0)$ of one of these critical submanifolds only a finite number of times, since $F$ is bounded below and decreases by a fixed amount each time.  Thus one can show that ${\mathcal M}$ is a disjoint union of the unstable manifolds of the nondegenerate critical submanifolds for $F$.

We remark that the gradient-like vector field ${\mathcal X}$ can be chosen to be tangent to each of the submanifolds $L^{2}_k(S^2,M)$, for $k \in {\mathbb N}$, $k \geq 2$, because we can choose local pseudogradients and the local representatives ${\mathcal X}_A$ near each critical submanfiold to satisfy this property.  Indeed, the only place in the argument where we need the $L^{2\alpha }_1$ topology is for convergence of trajectories to an $L^{2\alpha }_1$ critical point as $t \rightarrow \infty$.  Then elliptic bootstrapping shows that the critical points are in fact $C^\infty $, and convergence holds in all the Sobolev norms.  We formalize this in a remark:

\vskip .1 in
\noindent
{\bf Remark 4.1.}  We can choose our pseudogradient so that all critical points and all orbits between critical points lie within $L^{2}_k(S^2,M)$, for any $k \geq 2$.

\section{Morse index estimates}
\label{S:indexestimates}

We would like to show that our curvature condition implies that there are only finitely many orbits of nonconstant prime minimal two-spheres with Morse index less than a given bound.   The following estimates (analogous to the Morse-Schoenberg theorem for geodesics found in \cite{GKM}) provide a first step in this direction:

\vskip .1 in
\noindent
{\bf Theorem 5.1.} {\sl Suppose that the compact Riemannian manifold $M$ has positive half-isotropic sectional curvature.
\begin{enumerate}
\item Then any harmonic two-sphere has Morse index at least $n-2$.
\item There is a constant $c > 0$ such that if $f :S^2 \rightarrow M$ is any harmonic map, then
$$(\hbox{Morse index of $f$}) + 1 \geq c E(f).$$
\end{enumerate}}

\vskip .1 in
\noindent
{\bf Proof of 1:}  By a theorem of Koszul and Malgrange \cite{KoMa}, the bundle ${\bf E} = f^*TM \otimes {\mathbb C}$ has a canonical holomorphic structure, and the equation for harmonic maps implies that $\partial f/\partial z$ is a holomorphic section of this bundle.  Moreover, since the Riemannian metric is invariant under the Levi-Civita connection, it extends to a complex bilinear form
$$\langle \cdot , \cdot \rangle : {\bf E} \times {\bf E} \longrightarrow {\mathbb C},$$
which is holomorphic as a section of $\hbox{Hom}({\bf E} \otimes {\bf E},{\mathbb C})$.  Let ${\bf L}$ be the holomorphic subbundle of ${\bf E}$ generated by $\partial f/\partial z$ with orthogonal complement
$${\bf L}^\bot = \bigcup _{p \in S^2} \{ v \in {\bf E}_p : \langle v, \partial f/\partial z \rangle (p) = 0 \}.$$
Then ${\bf L} \subseteq {\bf L}^\bot$, so we can define the holomorphic normal bundle to be
\begin{equation} {\bf N} = {\bf L}^\bot / {\bf L},\label{E:normal}\end{equation}
and $\langle \cdot , \cdot \rangle$ descends to a holomorphic complex bilinear form on ${\bf N}$.  By a well-known theorem of Grothendieck \cite{Gk}, we can write ${\bf N}$ as a holomorphic direct sum of line bundles,
$${\bf N} = {\bf L}_1 \oplus {\bf L}_2 \oplus \cdots \oplus {\bf L}_{n-2}, \ \  \hbox{where} \ \ c_1({\bf L}_1)[S^2] \geq c_1({\bf L}_2)[S^2] \geq \cdots \geq c_1({\bf L}_{n-2})[S^2].$$
Moreover, the existence of a holomorphic complex bilinear form implies that
$$c_1({\bf L}_i)[S^2] = - c_1({\bf L}_{n-i-1})[S^2].$$

The second variation formula (\ref{E:secondvariation0}) for $E$ shows that
$$d^2E(f)(V,\bar V) \leq 2 \int _{S^2} \left[ \left\| \bar \partial V \right\|^2 - \kappa e(f) \langle V, \bar V \rangle \right] dA, \quad \hbox{for} \quad V \in \Gamma ({\bf N}),$$
where $\kappa $ is a positive lower bound on the half-isotropic curvature.  By the Riemann-Roch theorem (see the discussion in \S 4.1.1 of \cite{Mo2}) the space ${\mathcal O}({\bf L}_i)$ of holomorphic sections has dimension at least two when $c_1({\bf L}_i)[S^2]$ is positive and has dimension at least one when $c_1({\bf L}_i)[S^2] = 0$.  It follows that the total dimension of the space of holomorphic sections is at least $n-2$.  If $V$ is such a holomorphic section, then $\hbox{span}(V, \partial f/\partial z)$ is half-isotropic, so the assumption of positive half-isotropic curvature implies that the Morse index of $f$ is at least $n-2$, proving assertion 1.

\vskip .1 in
\noindent
{\bf Proof of 2:}  We need to sharpen the previous estimate and show that the Morse index of a harmonic two-sphere grows linearly with energy.  We can restrict to sections of one of the line bundle summands ${\bf L}_i$ of ${\bf N}$, obtaining
$$d^2E(f)(V,\bar V) \leq \int _{S^2} \left[ 2 \left\| \bar \partial V \right\|^2 - \phi \langle V, \bar V \rangle \right] dA, \quad \hbox{for} \quad V \in \Gamma ({\bf L}_i),$$
where $\phi = 2 \kappa e(f) > 0$.  To estimate the index of the restriction to $\Gamma ({\bf L}_i)$, we need to estimate the spectrum of the Schr\"odinger operator
$$2 \bar \partial ^* \circ \bar \partial - \phi : C^\infty ({\bf L}_i) \rightarrow C^\infty ({\bf L}_i),$$
where $\phi > 0$.  The eigenvalues of $2 \bar \partial ^* \circ \bar \partial - \phi$ are related to the eigenvalues of a scalar multiple of the operator $2 \bar \partial ^* \circ \bar \partial$ by the following reduction, which comes from Li and Yau \cite{LY}:  Since
\begin{equation} \frac{\int _{S^2} [ 2 |\bar \partial V |^2 - \phi |V| ^2] dA}{ \int _{S^2} |V| ^2 dA} = \frac{ \int_{S^2} \phi |V| ^2 dA}{\int _{S^2} |V|^2 dA} \left[ \frac{\int _{S^2} 2 |\bar \partial V |^2 dA}{\int_{S^2} \phi |V|^2 dA} - 1\right], \end{equation}
the index of the restriction of $d^2E(f)$ to $\Gamma ({\bf L}_i)$ is the maximal dimension of a linear subspace of $C^\infty ({\bf L}_i)$ on which the quotient
$$V \quad \mapsto \quad \frac{\int _{S^2} 2 |\bar \partial V |^2 dA}{\int_{S^2} \phi |V|^2 dA}$$
is strictly less than one.  We now divide into cases:

\vskip .1in
\noindent
{\bf $c_1({\bf L}_i) = 0$:}  In this case, we can choose the unitary structure on ${\bf L}_i$ so that
$$\frac{\int _{S^2} 2 |\bar \partial V |^2 dA}{\int_{S^2} \phi |V|^2 dA} = \frac{\int _{S^2} |dV |^2 dA}{\int_{S^2} \phi |V|^2 dA}.$$
But the right-hand side is just the Rayleigh quotient for the Laplace-Beltrami operator $\Delta _\phi$ for the Riemannian metric
$$\phi ds_0^2, \quad \hbox{where $ds_0^2$ is the base Riemannian metric on $S^2$.}$$
It follows from Korevaar's theorem \cite{K} that there is a constant $C > 0$ such that the $k$-th eigenvalue of $\Delta _\phi$ is less than or equal to one if
$$C\frac{k}{\int_{S^2} \phi dA} \leq 1, \quad \hbox{which occurs if} \quad k \leq \frac{1}{C}\int_{S^2} \phi dA.$$
Alternatively, we can say that the number of eigenvalues which is less than one is the greatest integer less than
$$C_1 \int_{S^2} \phi dA \quad \hbox{where} \quad C_1 = \frac{1}{C}.$$
This implies that the index of the Schr\"odinger operator $2 \bar \partial ^* \circ \bar \partial - \phi$ grows linearly with the integral of $\phi$, which gives the estimate of Theorem~5.1 in this case.

\vskip .1in
\noindent
{\bf $c_1({\bf L}_i) \neq 0$:}  (Since the sum of the Chern classes is zero, we need only treat the case in which $c_1({\bf L}_i) > 0$.)  To treat the case of Chern class $k > 0$, we use the fact that the spectrum of the Laplace operator is not affected much by what happens in the $\varepsilon $-ball around a single point, and note that if we take a single point out of the base of a $U(1)$-bundle, we lose the first Chern class, a fact that we can exploit.  Indeed, suppose that ${\bf E}_0$ is the trivial line bundle with trivial holomorphic structure.  We can let ${\bf E} = {\bf E}_0 \otimes \nu _p^k$, where $\nu _p$ is the point bundle at $p$ (as described for example in \S 7c of \cite {Gu}), the resulting bundle having a preferred holomorphic section $\phi _0$ which has divisor $(\phi _0) = k p$.  We can then define a map
\begin{equation} \iota : C^\infty(S^2 , {\bf E}_0) \longrightarrow C^\infty(S^2 , {\bf E}) \quad \hbox{by} \quad \iota (\phi ) = \phi \phi _0. \label{E:iota} \end{equation}
Finally, we let $A_0$ be a unitary connection on the $U(1)$-bundle ${\bf E}$ which has constant curvature and has $\phi _0$ as a holomorphic section.  We claim that if the $m$-th eigenvalue $\lambda _m^0$ of the operator $2 (\bar \partial ^* \circ \bar \partial)$ on ${\bf E}_0$ is $\leq C$, then so is the $m$-th eigenvalue $\lambda _m$ of the operator $2 \bar \partial_{A_0} ^* \circ \bar \partial_{A_0}$ on ${\bf E}$.

Indeed, if $\lambda _m^0 \leq C$, then there is an $m$-dimensional subspace ${\mathcal V}$ of $C^\infty(S^2 , {\bf E}_0)$ such that
$$\sigma \in {\mathcal V} \quad \Rightarrow \quad \frac{\int _{S^2} \langle \bar \partial \sigma, \bar \partial \sigma \rangle dA}{\int _{S^2} \left| \phi \right| ^2dA} \leq C.$$
We now look at the Rayleigh quotients of the corresponding elements $ \iota (\sigma )$ in $C^\infty(S^2 , {\bf E})$ and exploit the fact that
$$ \bar \partial_{A_0}(\sigma \sigma _0) = \bar \partial ( \sigma) \sigma _0.$$
This implies that the corresponding Rayleigh quotient in ${\bf E}$ is
$$\frac{\int _{S^2} \langle \bar \partial \sigma, \bar \partial \sigma \rangle |\sigma _0|^2 dA}{\int _{S^2} \left| \phi \right| ^2 |\sigma _0|^2dA}.$$
But this is just the Rayleigh quotient on the trivial bundle, with the unitary structure multiplied by $| \sigma _0 |^2$.  This would have the same spectrum as the standard unitary structure, except for the fact that $| \sigma _0 |^2$ vanishes at the point $p$.  Indeed, we can approximate $| \sigma _0 |^2$ by a nowhere zero function $\mu ^2$ in such a way that the support of $| \sigma _0 |^2 - \mu ^2$ is contained in an $\varepsilon $-ball around $p$.  Multiplication by a suitable cutoff function shows that this alters the spectrum by a small amount which goes to zero as $\varepsilon \rightarrow 0$, as we next describe.

\vskip .1 in
\noindent
{\bf Construction of a cutoff function:}  Following Choi and Schoen \cite{CS}, we define a map $\phi _\varepsilon: {\mathbb R} \rightarrow {\mathbb R}$ by
\begin{equation} \phi _\varepsilon (r) = \begin{cases} 0, & \hbox{if $r \leq \epsilon ^2$,} \cr 2 - (\log r)/(\log \epsilon ), & \hbox{if $\epsilon ^2 \leq r \leq \epsilon $,} \cr 1, & \hbox{if $\epsilon  \leq r$,} \end{cases} \label{E:phiepsilon}\end{equation}
so that
$$\frac{d\phi _\varepsilon}{dr} (r) = \begin{cases} 0, & \hbox{if $r \leq \epsilon ^2$,} \cr (-1)/(r\log \epsilon ), & \hbox{if $\epsilon ^2 \leq r \leq \epsilon $,} \cr 0, & \hbox{if $\epsilon  \leq r$,} \end{cases}$$
and
$$\int _0^{2\pi}\int _0^\epsilon \left( \frac{d\phi _\varepsilon}{dr} (r) \right) ^2 rdrd\theta = \int _{\epsilon ^2}^\epsilon \frac{2\pi}{r (\log \epsilon )^2}dr = - \frac{2 \pi}{\log \epsilon }.$$
Thus if we define $\psi _\varepsilon : \Sigma \rightarrow {\mathbb R}$ so that it is one outside an $\epsilon$-neightborhood of the $p$ and in terms of geodesic polar coordinates $(r, \theta )$ about $p$ (for the constant curvature metric) satisfies the condition $\psi _\varepsilon = \phi _\varepsilon \circ r$, then
$$\int _{S^2} |d \psi _\varepsilon |^2 dA \leq \frac{-C}{\log \epsilon }, \qquad \hbox{where $C$ is a positive constant,}$$
and $\psi _\varepsilon$ is our desired cutoff function which vanishes in the $\varepsilon ^2$ neighborhood of the point $p$, and yet alters Rayleigh quotients by arbitrarily little.

\vskip .1 in
\noindent
{\bf Remark 5.2.}  In Part 1 of this theorem that we use the hypothesis of positive half-isotropic sectional curvatures; for manifolds of positive isotropic curvature, we would get a weaker estimate on the Morse index.  It follows from the proof of the Sphere Theorem of \cite{MM} that any compact simply connected Riemannian manifold of dimension $n \geq 4$ always has nonconstant minimal two-spheres of Morse index $\leq n-2$, so such manifolds with positive half-isotropic curvature always have minimal two-spheres of Morse index exactly $n-2$.

For Part 2 of this theorem we only need positive curvature for the half-isotropic two-plane generated by $\partial f/\partial z$ and one of the line bundle summands in the holomorphic direct sum decomposition of the normal bundle.  Thus, for example, Part 2 holds for four-dimensional manifolds which satisfy the weaker condition of positive isotropic curvature.

\section{Bubble trees and the Nonbubbling Theorem}
\label{S:nonbubbling}

We are interested in critical points for the usual energy $E$, not the $\alpha $-energy $E_\alpha $ or the Morse approximation $F_\alpha $ to $E_\alpha $ defined by (\ref{E:definitionF}).  Thus we need to study what happens to the critical points of $E_\alpha $ or $F_\alpha $ as the perturbations are turned off.

This is a special case of a more general question asked by Sacks and Uhlenbeck in \cite{SU1}.  If $\Sigma $ is a compact Riemann surface with conformal structure $\omega $, the $\omega $-energy is the function
$$E_\omega : \hbox{Map}(\Sigma ,M) \longrightarrow M \quad \hbox{defined by} \quad E_\omega (f) = \frac{1}{2} \int _{\Sigma } |df|^2 dA,$$
where $|df|$ and $dA$ are defined in terms of a Riemannian metric in the conformal equivalence class determined by $\omega $, and it has a perturbation
$$E_{\alpha ,\omega} : \hbox{Map}(\Sigma ,M) \longrightarrow M \quad \hbox{defined by} \quad E_{\alpha ,\omega} (f) = \frac{1}{2} \int _\Sigma [(1 + |df|^2)^\alpha - 1]dA.$$
In the perturbation we lose conformal invariance; we often choose the metric on $\Sigma $ to have constant curvature and to satisfy a normalization condition, such as having total area one.

Suppose that for $m \in {\mathbb N}$, $f_m$ is a minimax critical point for $E_{\alpha _m,\omega}$ with energy bounded by $E_0 > 0$ such as might be constructed for a given homology class when working out a Morse theory for $E_{\alpha _m,\omega}$.  Does the sequence $\{ f_m : m \in {\mathbb N} \}$ have a subsequence which has a harmonic map limit as $m \rightarrow \infty$?  Sacks and Uhlenbeck show that there is indeed such a subsequence which converges in a weak sense to a base harmonic map $f_0 : \Sigma \rightarrow M$ together with a collection of harmonic two-spheres which bubble off in the limit.  Later, Parker and Wolfson introduced the term of \lq\lq bubble tree" to describe the limit, and the nature of the convergence was further studied by Parker \cite{Par}, and Chen and Tian \cite{CT}.  The theory of such bubble trees is summarized in \S 4.6 and \S 4.9 of \cite{Mo2}. We are interested here in the case in which $\Sigma = S^2$.
 
\vskip .1in
\noindent
{\bf Definition.}  A {\em bubble tree} based upon the Riemann surface $S^2$ is a collection of maps indexed by the vertices and edges of a rooted tree $T$.  To each element $v$ in the vertex set $V$, there corresponds a minimal two-sphere
$$g_v : S^2 \rightarrow M,$$
which may in some cases reduce to a constant.  The root is denoted by $0$ and corresponds to a {\em base map} which is also a minimal two-sphere
$$g_0 : S^2 \rightarrow M,$$
while the other minimal two-spheres are called {\em bubbles}.
On the other hand, the edges correspond to smooth geodesics called {\em necks} 
$$\gamma _v :[0,1] \rightarrow M$$
which join the parent $p(v)$ of a vertex $v$ to $v$ itself (with $\gamma _v$ only defined when $v \neq 0$), so that
$$\gamma _v(0) \in \hbox{Image}(f_{p(v)}), \quad \gamma _v(1) \in \hbox{Image}(f_{v}).$$
We also allow some of the necks $\gamma _v$ to be constant maps to a single point of $M$.  Finally, we require that $f_v$ be nonconstant when $v$ is not the parent of another vertex.

\vskip .1 in
\noindent
In what sense does a subsequence of the sequence $\{ f_m : m \in {\mathbb N} \}$ of critical points or
$$E_{\alpha _m} : \hbox{Map}(S^2 ,M) \longrightarrow {\mathbb R}, \quad \alpha _m \rightarrow 1$$
converge to a bubble tree?  We can give a schematic description with the above references providing details.  We replace the sequence by its subsequence and for each $m \in {\mathbb N}$ we construct a family of disjoint disks $\{ D_{m,v} : v \in V \}$ within $S^2$ on which energy density is concentrating.  When chosen appropriately we can choose $f_m|D_{m,v}$ so that it rescales to a mapping
$$\tilde f_{m,v} : \{ z \in {\mathbb C} : |z| \leq m \} \longrightarrow M$$
converges to a harmonic parametrization of
$$f_v|(S^2 - \{ \infty , p_1, \ldots p_k \}),$$
where $f_v(p_1), \ldots f_v(p_k)$ are endpoints of geodesics corresponding to edges linking $v$ to further vertices $v_1, \ldots v_k$.  The convergence of the rescaling $\tilde f_{m,v}$ of $ f_m|D_{m,v}$ to $f_v$ is uniform on compact subsets of $S^2 - \{ \infty , p_1, \ldots p_k \}$.  One shows that the restriction of $f_m$ to $$S^m - \bigcup_{v\in V} D_{m,v}$$
converges in a suitable sense to the geodesics corresponding to edges in the bubble tree.

The curvature assumptions positive complex sectional curvatures or positive half-isotropic curvatures imply that Ricci curvature is positive, hence the fundamental group of $M$ is finite by the Theorem of Myers.  Therefore, Theorem~4.9.2 of \cite{Mo2} implies that necks are geodesics of finite length and zero energy, so no energy is \lq\lq lost in the necks."  In particular, after passing to a subsequence, $E_{\alpha _m}$ converges to the sum of the areas of the nonconstant minimal two-spheres $f_v$ obtained in the limit.

We emphasize that the only elements of the bubble tree which can carry nonzero energy are the nonconstant minimal two spheres $f_v$.  Moreover, if we assume that the real sectional curvature of $M$ satisfy the inequality $K_r(\sigma ) \leq 1$, then an upper bound on total energy gives a limit on the number of nonconstant bubbles obtained.  (We can always rescale the metric on a compact manifold so that it satisfies this condition.)  Indeed, it follows from the Gauss equation that if $f : S^2 \rightarrow M$ is any minimal two-sphere, the induced Gaussian curvature $K_f(p)$ at $p \in S^2$ satisfies
\begin{equation} K_f(p) \leq K_r((f_p)_*T_pS^2) \label{E:Gauss}\end{equation}
so if $K_r(\sigma ) \leq 1$, the Gauss-Bonnet formula implies
$$E(f) = \int _{S^2} dA_f \geq \int _{S^2} K_fdA_f = 4 \pi.$$
Thus each nonconstant minimal two-sphere $f_v$ has area and energy at least $4\pi$, and if the total energy is bounded above by $a > 0$, the number of nonconstant minimal two-spheres in any bubble tree is $\leq a/4\pi$.

\vskip .1 in
\noindent
{\bf Lower Semicontinuity Theorem~6.1.} {\sl Suppose that $\{ f_m : m \in {\mathbb N} \}$ is a sequence of critical points for $E_{\alpha _m}$ which converge to a bubble tree in the sense described above.  Then the Morse indices of the nonconstant minimal two-spheres in the bubble tree satisfy
\begin{equation} \sum \{ \hbox{Morse index of $f_v$} : \hbox{$f_v$ is nonconstant} ) \leq \hbox{Morse index of $f_m$}, \label{E:lowersemicontinuity}\end{equation}
when $m$ is sufficiently large.}

\vskip .1in
\noindent
Proof:  Our strategy is to show that if $d^2E_\omega (f_v)$ is negative definite on a linear space ${\mathcal V}$ of dimension $q$, then its approximation $d^2E_{\alpha _m,\omega} (f_m)$ is negative definite on a linear space of sections with compact support in $D_{m,v}$ of dimension $q$.  This proves our claim since the $D_{m,v}$'s are disjoint.  

To determine the Morse index of $f_m$ we use the second variation formula for $\alpha $-harmonic maps $f$:
\begin{multline*} d^2E_{\alpha ,\omega} (f)(V,W) = \alpha \int_{S^2} (1 + |df|^2)^{\alpha -1} [\langle \nabla V ,  \nabla W \rangle - \langle {\cal K}(V), W \rangle ] dA \\ + 2\alpha (\alpha - 1) \int_{S^2} (1 + |df|^2)^{\alpha -2} \langle df, \nabla V \rangle \langle df, \nabla W \rangle dA.\end{multline*}
Note that by inequality (\ref{E:quotientestimate}),
$$2\alpha (\alpha - 1) \int_{S^2} (1 + |df|^2)^{\alpha -2} \langle df, \nabla V \rangle \langle df, \nabla W \rangle dA \leq \varepsilon \alpha \int_{S^2} (1 + |df|^2)^{\alpha -1} \langle \nabla V ,  \nabla W \rangle dA,$$
when $\alpha $ is sufficiently close to one.  After rescaling, the first term approaches $d^2E_\omega (f)(V,W)$ when the supports of $V$ and $W$ do not contain any of the bubble points $\{ p_0 = \infty, p_1, \ldots , p_k \}$.

To show that the bubble points don't interfere, we make use of a cutoff function near the bubble points, defined utilizing the function $\phi _\varepsilon$ of (\ref{E:phiepsilon}).  We define $\psi _i : S^2 \rightarrow {\mathbb R}$ so that it is one outside an $\epsilon$-neightborhood of the bubble point $p_i$ and in terms of polar coordinates $(r_i, \theta _i)$ about the bubble point $p_i$ satisfies the condition $\psi _i = \phi _\varepsilon \circ r_i$, so that
$$\int _{S^2} |d \psi _i |^2 dA \leq \frac{-C}{\log \epsilon }, \qquad \hbox{where $C$ is a positive constant,}$$
and check that a similar estimate holds for $\psi = \psi _0 \psi _1 \cdots \psi _k$, a cutoff function which vanishes at every bubble point.  If $\varepsilon > 0$ is chosen sufficiently small, then
$$V \in {\mathcal V} \quad \Rightarrow \quad d^2E(f_v)(\psi V, \psi V) < 0 \quad \Rightarrow \quad d^2E_\alpha (f_m)(\psi V, \psi V) < 0.$$
This shows that $d^2E_\alpha (f_m)$ is negative definite on a space of dimension $q$ with support lying in $D_{m,v}$.  QED

\vskip .1in
\noindent
The Nonbubbling Theorem~1.4 follows from Theorem~6.1 and Theorem 5.1.1, the latter theorem stating that any harmonic two-sphere in a manifold with positive half-isotropic sectional curvatures has Morse index at least $n - 2$.  Thus if a bubble tree contains $k$ bubbles of nonzero area, it must have Morse index at least $k(n-2)$.

%Sept25 replace statement of compactness theorem and insert prime in the proof:

%\vskip .1 in
%\noindent
%{\bf Compactness Theorem~6.2.} {\sl  Suppose that $M$ has a generic Riemannian metric with positive half-isotropic sectional curvatures.  Then there are only finitely many ${\mathcal G}$-orbits of nonconstant prime minimal two-spheres $f : S^2 \rightarrow M$ of Morse index $\leq 2n-5$, where ${\mathcal G} = PSL(2,{\mathbb C}) \cup R \cdot PSL(2,{\mathbb C})$.}

\vskip .1 in
\noindent
{\bf Compactness Theorem~6.2.} {\sl  Suppose that $M$ has a generic Riemannian metric with positive half-isotropic sectional curvatures, and $C$ is a closed subset of $C^k(S^2,M)$ which contains no nontrivial branched covers.  Then there are only finitely many ${\mathcal G}$-orbits of nonconstant prime minimal two-spheres $f : S^2 \rightarrow M$ of Morse index $\leq 2n-5$ in $C$, where ${\mathcal G} = PSL(2,{\mathbb C}) \cup R \cdot PSL(2,{\mathbb C})$.}

\vskip .1 in
\noindent
Proof:  Suppose the contrary, that there is a sequence $\{ f_m \}$ of minimal two-spheres in $C$ of Morse index  $\leq 2n-5$ lying on distinct ${\mathcal G}$-orbits.  By Theorem~5.1.2, these must all have energy $< c$, where $c$ depends on the maximum possible Morse index $2n-5$.  We can then use the Perturbation Theorem~2.2 to construct a corresponding sequence $\{ f'_m \}$ of $\alpha _m$-energy critical points for $E_{\alpha _m}$, with bounded $\alpha _m$-energy and $\alpha _m \rightarrow 1$.  A subsequence of these must converge to a harmonic map $f_\infty : S^2 \rightarrow M$ without bubbling, because bubbling would contradict the Nonbubbling Theorem~1.4.  But then we have a sequence of harmonic maps on distinct ${\mathcal G}$-orbits which converges to $f_\infty $ contradicting the nondegeneracy of the limit.

\section{Branched covers}
\label{S:branchedcovers}

Suppose that $h : S^2 \rightarrow M$ is an immersed prime minimal two-sphere.  If $g : S^2 \rightarrow S^2$ is a holomorphic map (also known as a meromorphic function on $S^2$), then the composition $f = h \circ g: S^2 \rightarrow M$ is a conformal harmonic map, hence a minimal surface in its own right, which we call a {\em branched cover\/} of $h$.  To study these, we consider the space 
$$\hbox{Hol}_d(S^2 ,S^2) = \{ g: S^2 \rightarrow S^2 : \hbox{ $g$ is holomorphic and has degree $d$ } \},$$
for each integer $d \geq 2$.  For each integer $d \geq 1$, we define a map
$$\Phi _d : \hbox{Map}(S^2 ,M) \times \hbox{Hol}_d(S^2 ,S^2) \rightarrow \hbox{Map}(S^2 ,M) \quad \hbox{by} \quad \Phi _d(h,g) = h \circ g ,$$
a map which will be $C^k$ when the mapping spaces are completed with respect to suitable topologies.  If $h$ is an immersed prime minimal surface, the image of $\Phi _d$, the space of all branched covers of degree $d$, is a critical submanifold
$${\mathcal M}_d(h) \subseteq L^2_k(S^2,M).$$
However, Example~2.8.2 from \cite{Mo2} shows that the Morse index can vary from point to point of ${\mathcal M}_d(h)$, so the space ${\mathcal M}_d(h)$ of branched covers is not necessarily a nondegenerate submanifold in the sense of Bott.

Our approach is to utilize the pinching condition to show that branched covers of degree three or larger have Morse index at least $2n-4$, and then study an equivariant perturbation of the $\alpha$-energy near two-fold branched covers.  The first step is given by the following estimate:

\vskip .1 in
\noindent
{\bf Theorem 7.1.} {\sl Suppose $M$ is an $n$-dimensional Riemannian manifold which has half-isotropic and real sectional curvatures which satisfy the curvature condition
\begin{equation} K_i(\sigma)>\frac{1}{d} K_r(\widehat\sigma) > 0, \label{E:curvd}\end{equation}
where $\widehat\sigma$ is a real two-plane associated to $\sigma $.  Then the Morse index of any branched cover of degree $\geq d$ of a minimal two-sphere $f: S^2 \rightarrow M$ must be at least $2n-4$.}

\vskip .1 in
\noindent
Our pinching condition (\ref{E:1/3bis}) gives (\ref{E:curvd}) when $d = 3$, so that under assumption (\ref{E:1/3bis}) all branched covers of degree at least three have Morse index at least $2n-4$.

\vskip .1 in
\noindent
Proof of Theorem 7.1:  As in the proof of Theorem~5.1, we use the normal bundle (\ref{E:normal}) and its holomorphic decomposition into line bundles,
$${\bf N} = {\bf L}_1 \oplus {\bf L}_2 \oplus \cdots \oplus {\bf L}_{n-2}, \ \  \hbox{where} \ \ c_1({\bf L}_1)[S^2] \geq c_1({\bf L}_2)[S^2] \geq \cdots \geq c_1({\bf L}_{n-2})[S^2].$$
If $V = X - i Y$ is a holomorphic section of a line bundle summand ${\bf L}_{i}$ with $c_1({\bf L}_i)[S^2] > 0$, then $V$ must vanish at some point, and since $\langle V,V \rangle$ is a bounded holomorphic function and hence constant, $\langle V,V \rangle \equiv 0$.  This implies that the real and imaginary components $X$ and $-Y$ must have the same length and be perpendicular to each other.  If $E_i$ denotes the rank two real subbundle of $f^*TM$ generated by real and imaginary parts of sections of ${\bf L}_{i}$, we can define an almost complex structure $J$ on $E_i$ by
$$J(X) = Y, \quad \hbox{whenever} \quad V = X - i Y \in \Gamma({\bf L}_{i}).$$
Then the components $X$ and $Y$ of a holomorphic section $V = X - i Y$ of ${\bf L}_{i}$ generate a two-dimensional subspace of sections of $E_i$, complex linear with respect to $J$, on which the Morse index is negative.  Since 
$$c_1({\bf L}_i)[S^2] = k > 0 \quad \Rightarrow \quad \dim {\mathcal O}({\bf L}_i) = k + 1,$$
the corresponding real and imaginary parts form a linear subspace of sections of real dimension $2k+2$ on which the index form is negative.  The upshot is that each summand with positive Chern class must contribute at least four to the Morse index.

If there are $m$ summands with positive first Chern class (where $m \leq (n-2)/2$), there will be $n-2-2m$ summands with zero Chern class.  Thus if we can show that each zero Chern class summand contributes at least two to the Morse index, the total Morse index will be at least
$$4m + 2 \times (n-2-2m) = 2n-4,$$
establishing Theorem~7.1.

Thus we need only show that our curvature hypothesis (\ref{E:1/3bis}) implies that trivial line bundle summands must each contribute at least two to the Morse index.  The holomorphic section always gives a contribution of one to the index, and to obtain a second contribution of one, we apply the second variation formula (\ref{E:secondvariation0}) to a section corresponding to the first nonzero eigenvalue, obtaining
\begin{multline} d^2E(f)(V,\bar V) = 2 \int _{S^2} \left[ \left\| \bar \partial V \right\|^2 - K_i(\sigma) e(f) \left\| V \right\|^2 \right] dA \\ = \int _{S^2} \left[ 2 \left\| \bar \partial V \right\|^2 - 2  K_i(\sigma) \left\| V \right\|^2 \right] d\tilde A, \label{E:indexestimate7} \end{multline}
where $d \tilde A = e(f) dA$ is the area form for the induced metric on $S^2$, and $\sigma$ is the span of $V$ and $\frac{\partial f}{\partial z}$.  By Lemma~5.2, we can assume that the connection is flat and hence (\ref{E:indexestimate7}) simplifies to
$$d^2E(f)(V,\bar V) = \int _{S^2} \left[ \left\| dV \right\|^2 - 2  K_i(\sigma)\left\| V \right\|^2 \right] d\tilde A.$$
Substituting the curvature condition (\ref{E:curvd}), we have
\begin{multline} d^2E(f)(V,\bar V) < \int _{S^2} \left[ \left\| dV \right\|^2 - \frac{2}{d}  K_r(\widehat \sigma)\left\| V \right\|^2 \right] d\tilde A,\\
=\left(\frac{\int_{S^2}\left\| dV \right\|^2  d\tilde A}{\int_{S^2} K_r(\widehat \sigma)\left\| V \right\|^2  d\tilde A}-\frac{2}{d}\right)\int_{S^2} K_r(\widehat \sigma)\left\| V \right\|^2  d\tilde A, \\ \leq\left(\lambda_1-\frac{2}{d}\right)\int_{S^2} K_r(\widehat \sigma)\left\| V \right\|^2  d\tilde A,\end{multline} 
where $\lambda _1$ denotes the first nonzero eigenvalue of the Laplace-Beltrami operator for the metric $K_r(\widehat\sigma)f^* g$ on $S^2$ which has area form $K_r(\widehat \sigma)  d\tilde A$.  Applying Hersch's Theorem (Theorem~1 of Yang and Yau \cite{YY} or Theorem 7.6 of Colding and Minnicozzi \cite{CM}), we find that $\lambda_1$ satisfies
$$\lambda _1 \leq \frac{8\pi}{\int_{S^2} K_r(\widehat \sigma)  d\tilde A}\leq \frac{8\pi}{\int_{S^2} K_f  d\tilde A}\leq\frac{8\pi}{4\pi d} = \frac{2}{d},$$
for a $d$-fold branched cover, where we have used that $\widehat\sigma=f_* TS^2$ and, as in (\ref{E:Gauss}), $K_f$ is the Gauss curvature of the induced metric which satisfies, $K_f\leq K_r(f_* TS^2)$.  Thus we find that there are at least two variation fields for each trivial summand of the normal bundle which decrease area, finishing the proof.

\vskip .1 in
\noindent
In the special case in which $M = S^n$ has the round metric of constant curvature one, Theorem~7.1 implies:

\vskip .1 in
\noindent
{\bf Corollary 7.2.} {\sl If $M = S^n$ has the round metric of constant curvature one and $h : S^2 \rightarrow S^n$ is a totally geodesic imbedding, then the Morse index of any double branched cover $f \in {\mathcal M}_2(h)$ is $\geq 2n-4$.}

\vskip .1 in
\noindent
Proof:  In this case, (\ref{E:curvd}) holds for all $d \geq 2$.

\vskip .1 in
\noindent
For metrics of nonconstant curvature, our curvature condition is not strong enough to ensure Morse index at least $2n-4$ for two-fold branched covers, and we will use nontrivial isotropy of the group action to deal with these two-fold branched covers.  The analysis of two-fold branched covers simplifies because any element of $g \in \hbox{Hol}_2(S^2 ,S^2)$ can be related to the standard two-fold cover $\pi _2$,
\begin{equation} g = S^{-1} \circ \pi _2 \circ T, \quad \hbox{where} \quad S,T \in PSL(2,{\mathbb C}) \quad \hbox{and} \quad \pi _2(z) = z^2. \label{E:doublecover}\end{equation}
Recalling the discussion on center of mass from \S \ref{S:perturbations}, we reduce the structure group on ${\mathcal M}_2(h)$ from $\widehat {\mathcal G}$ to $\widehat G = O(3)$ by restricting to the submanifold
$${\mathcal M}_2(h)_0 = \left\{ f \in {\mathcal M}_2(h) : \int _{S^2} {\bf X} \psi _1(|df|^2) dA = {\bf 0} \right\}$$
of branched covers with center of mass zero, the action of $\widehat G = O(3)$ on ${\mathcal M}_2(h)_0$ preserving the center of mass zero condition.  It follows from the explicit form (\ref{E:doublecover}) of double branched covers that any element of ${\mathcal M}_2(h)_0$ must have nontrivial isotropy for the action of $G = SO(3)$, the isotropy group being ${\mathbb Z}_2$, generated by a rotation through angle $\pi $.  Milnor \cite{Mil3} describes the quotient of ${\mathcal M}_2(h)_0$ by the $SO(3)$ action as a four-manifold which is the total space of a two-plane bundle over ${\mathbb R}P^2$.

Although ${\mathcal M}_2(h)_0$ is not compact, we can make it compact by removing the branched covers whose two branch points are too close.  Given $\varepsilon _0 >0$, we set
$${\mathcal M}_2(h)_0^{\varepsilon _0} = \{ f \in {\mathcal M}_2(h)_0: \hbox{ the two branch points have distance $\geq \varepsilon _0$ } \},$$
the distance being measured in the standard metric of constant curvature and total area one on $S^2$.  We can describe the topological type of this space as follows:  For the given fixed prime minimal immersion $h : S^2 \rightarrow M$ we can specify a $\widehat G$-orbit of branched covers of $h$ by locating the images of its two branch points.  Let $n \in S^2$ be the first image, and regard
$$S^2 - (\hbox{ $\varepsilon _0$-neighborhood of $n$ })$$
as a two-disk $D_n$.  The second image is then a point in $D_n$ and the image of the branch locus is specified by a point in a $D^2$-bundle over $S^2$.  But we can interchange the two branch points and a quotient by the resulting ${\mathbb Z}_2$ action gives us the two-disk bundle $N^4$ over ${\mathbb R}P^2$.  Then ${\mathcal M}_2(h)_0^{\varepsilon _0}$ itself is the total space of a fibration
\begin{equation} \pi : {\mathcal M}_2(h)_0^{\varepsilon _0} \longrightarrow N^4 \label{E:projetoN}\end{equation}
the fiber being a quotient space $O(3)/{\mathbb Z}_2$, where ${\mathbb Z}_2$ is a subgroup but not normal.  Each element of ${\mathcal M}_2(h)_0^{\varepsilon _0}$ has ${\mathbb Z}_2$ isotropy under the $SO(3)$ action, while the larger $O(3)$ action interchanges the two components of ${\mathcal M}_2(h)_0^{\varepsilon _0}$.

We next study the $\alpha $-energy critical points which are close to ${\mathcal M}_2(h)_0$, for a given choice of $h : S^2 \rightarrow M$.  We can use the $L^2_k$ topology where $k$ is large, so our function spaces are Hilbert manifolds, as described in Lang \cite{L}.  The Hilbert space inner product on $L^2_k(S^2, {\mathbb R}^N)$ makes it into a flat infinite-dimensional Riemannian manifold as described in Chapter VII of \cite{L}, and the isometric imbedding $M \subseteq {\mathbb R}^N$ induces an imbedding $L^2_k(S^2,M) \subseteq L^2_k(S^2, {\mathbb R}^N)$, which induces a Riemannian metric on $L^2_k(S^2,M)$ for which $\widehat G$ acts as isometries; this in turn induces a Riemannian metric on the finite-dimensional manifold ${\mathcal M}_2(h)_0$, a space on which $\widehat G$ acts with ${\mathbb Z}_2$ isotropy.  We now apply the tubular neighborhood theorem: for a suitable $\widehat G$-invariant open neighborhood $\nu {\mathcal M}_2(h)_0$ of the zero section in the normal bundle to ${\mathcal M}_2(h)_0$ within $L^2_k(S^2,M)$, we can define a smooth $\widehat G$-equivariant diffeomorphism
$$\phi : \nu {\mathcal M}_2(h)_0\longrightarrow U,$$
where $U$ is some $\widehat G$-invariant open tubular neighborhood of ${\mathcal M}_2(h)_0$ in $L^2_k(S^2,M)$.  The function $E_\alpha $ restricts to a $\widehat G$-invariant function on $U$.

If $f_0 \in {\mathcal M}_2(h)_0^{\varepsilon _0}$, we let
$$N(f_0) = \hbox{fiber of $\nu {\mathcal M}_2(h)_0^{\varepsilon _0}$ over $f_0$,}$$
a submanifold of finite codimension within $L^2_k(S^2,M)$, and we restrict the function $E_\alpha $ to this submanifold, $E_\alpha : N(f_0) \rightarrow {\mathbb R}$.  By an immediate extension of the argument for Theorem 2.2 we can prove:

\vskip .1 in
\noindent
{\bf Proposition 7.3.} {\sl There exists $\alpha _0 \in (1,\infty)$ and a smooth map
\begin{multline} \eta : {\mathcal M}_2(h)_0^{\varepsilon _0} \times [1,\alpha _0) \\ \longrightarrow \left\{ (f, \alpha ) \in L^{2}_{k}(S^2 ,M) \times [1,\alpha _0) : \int _{S^2} {\bf X} \psi _\alpha (|df|^2) dA = {\bf 0} \right\} \label{E:defofeta}\end{multline}
which satisfies the following conditions:
\begin{enumerate}
\item $\eta $ is $O(3)$-equivariant,
\item $\eta _\alpha (f) = \eta (f, \alpha )$ is the unique critical point for the restriction of $E_\alpha $ to $N(f_0)$ for $f_0 \in {\mathcal M}_2(h)_0^{\varepsilon _0}$ and $\alpha \in [1,\alpha _0)$, and this critical point is nondegenerate on $N(f_0)$,
\item there is an $\varepsilon _1 > 0$ such that if $f$ is a critical point for $E_\alpha $ such that $\| f - f_0 \|_{C^2} < \varepsilon _1$ for some $f_0 \in {\mathcal M}_2(h)_0^{\varepsilon _0 + \varepsilon_1}$, $f$ lies in the image of $\eta $, and
\item the isotropy group under the $SO(3)$ action of every element in the image of $\eta $ is ${\mathbb Z}_2$, a rotation through angle $\pi $ about some axis. 
\end{enumerate}}

\vskip .1 in
\noindent
For the third condition on $\eta $, note that if $f$ is an element of $N(f_0)$ and a critical point for $E_\alpha $, it must also be critical point for the restriction of $E_\alpha $ to $N(f_0)$.  This proposition shows that critical points of $E_\alpha $ which are sufficiently close to the submanifold of two-fold branched covers must necessarily be invariant under a rotation in $O(3)$ of 180 degrees.

But we would like those critical points to lie on nondegenerate critical submanifolds, each an orbit for the $\widehat G$-action.  We can arrange this by first constructing a suitable smooth function $h : N^4 \rightarrow {\mathbb R}$, then letting $\tilde h = h \circ \pi$, where $\pi $ is the projection of (\ref{E:projetoN}).  Since $\tilde h$ is constant on each fiber of $\pi $ it is invariant under the action of $\widehat G$.  Finally, we let $\psi _\alpha: U \rightarrow [0,1]$ be a smooth $\widehat G$-invariant function which is identically one on the image of the function $\eta _\alpha = \eta (\cdot ,\alpha )$ defined by (\ref{E:defofeta}), and let
\begin{equation}F_h = E_\alpha + \psi _\alpha \tilde h, \label{E:FsubH}\end{equation}
For small $h$ this is a $\widehat G$-equivariant perturbation of $E_\alpha $ every critical point of which has ${\mathbb Z}_2$ isotropy under the $\widehat G$ action.

We claim that for generic choice of $h$, the perturbed function $F_h$ will have all critical points lying on isolated $\widehat G$ orbits, each nondegenerate in the sense of Bott.  Indeed, the only critical points for $F_h$ which lie in $\phi (N(f_0))$ for some $f_0 \in {\mathcal M}_2(h)_0^{\varepsilon _0}$ must be critical points for the restriction of $E_\alpha $ to $\phi (N(f_0))$ and must therefore lie in the image of
$$\eta _\alpha :{\mathcal M}_2(h)_0^{\varepsilon _0} \longrightarrow U.$$
On this image $\psi_\alpha \equiv 1$, and
$$F_h \equiv E_\alpha + \tilde h,$$
is a small perturbation of $E_\alpha $ when $h$ is small.  At a critical point $\eta _\alpha (f_0)$ for $F_h$, the tangent space divides into a direct sum
$$T_{\eta _\alpha (f_0)}L^{2}_{k}(S^2 ,M) \cong V \oplus H \oplus N, \quad \hbox{where} \quad H \oplus N \cong T_{f_0}{\mathcal M}_2(h)_0,$$
and we arrange that $V$ is tangent to the fiber of $\nu {\mathcal M}_2(h)_0$ (the \lq\lq vertical space"), $H$ projects via (\ref{E:projetoN}) to the tangent space to $N^4$ (the \lq\lq horizontal space"), and that $N$ is tangent to the orbit of the $\widehat G$ action.  At a critical point, $DE_\alpha = 0$ and $d\tilde h = 0$.  Invariance under $\widehat G$ implies that $N$ is a null-space for both $D^2E_\alpha$ and $d^2\tilde h$.  Moreover, since $h$ is constant on each fiber of the normal bundle,
$$\nu {\mathcal M}_2(h)_0 \longrightarrow {\mathcal M}_2(h)_0,$$
the space $V$ must lie in the null space of $d^2\tilde h$.  By part 2 of Proposition 7.3, $D^2E_\alpha$ is nondegenerate on $V \times V$, and a generic choice of $\tilde h$ will make
$$d^2F_h = D^2E_\alpha + d^2\tilde h \quad \hbox{nondegenerate on} \quad (V \oplus H) \times (V \oplus H),$$
establishing the claim.  We can summarize what we have proven as follows.

\vskip .1 in
\noindent
{\bf Proposition 7.4.} {\sl Suppose that $h : S^2 \rightarrow M$ is an imbedded minimal two-sphere in the compact Riemannian manifold $M$ and that $U$ is a small open neighborhood of the space ${\mathcal M}_2(h)_0^{\varepsilon _0 + \varepsilon_1}$ of two-fold branched covers of $h$, where $\varepsilon _0$ and $\varepsilon _1$ can be taken arbitrarily small.  Then there is a perturbation $\tilde E_\alpha $ of $E_\alpha$ which agrees with $E_\alpha $ outside $U$, is $C^2$ close to $E_\alpha $, and satisfies the following condition:  All critical points for $\tilde E_\alpha $ lie on nondegenerate critical submanifolds, each an orbit for the $\widehat G$ action.  Moreover, each critical point has isotropy group ${\mathbb Z}_2$, generated by some rotation through $\pi $.}

\vskip .1 in
\noindent
We have left out the points in ${\mathcal M}_2(h)_0 - {\mathcal M}_2(h)_0^{\varepsilon _0}$, but when the half-isotropic curvature is positive, these must have Morse index $\geq 2n-4$ when $\varepsilon _0$ is sufficiently small, by the proof of the Nonbubbling Theorem~1.4.

\section{The topology of the mapping space}
\label{S:equivariant}

The strategy for proving the Main Theorem~1.3 is to first use equivariant Morse theory to calculate the low degree equivariant cohomology of the mapping space $\hbox{Map}(S^2 ,M)$ (carried out in this section), and then use equivariant Morse theory in the other direction (carried out in \S~\ref{S:existence}) to prove existence of minimal two-spheres in manifolds which satisfy our curvature hypotheses.

Since the $\alpha $-energy $E_\alpha : \hbox{Map}(S^2 ,M) \rightarrow {\mathbb R}$ is invariant under an action of the group $O(3)$ by isometries on the domain, the nonconstant critical points for $E_\alpha $ lie on $O(3)$-orbits of critical points, which can be measured by $O(3)$-equivariant Morse theory.  We now review the essential features of equivariant cohomology.

Any compact Lie group $G$ possesses a classifying space $BG$ together with a universal principal $G$-bundle with contractible total space $EG$ and right action $\alpha _E : EG \times G \rightarrow EG$.  If $X$ is a CW complex and $\alpha _X: X \times G \rightarrow X$ is a continuous right action, the diagonal action on the product $X \times EG$,
$$(\sigma, (x, e)) \mapsto ( \alpha _X(x,\sigma), \alpha _E(e,\sigma )),$$
is free, and we can therefore divide by this action to obtain the topological space $X_G$.  To be more precise, we set
$$X_G = (X \times EG)/G,$$
the set of equivalence classes in $X \times EG$ under the equivalence relation
$$(\alpha _X (x , \sigma  ),  \alpha _E(e, \sigma )) \sim (x, e), \qquad \hbox{for $x \in X$, $e \in EG$ and $\sigma \in G$.}$$
We then define the {\em equivariant cohomology\/} of $X$ with coefficients in the coefficient ring $R$ by
$$H^*_G(X;R) = H^*(X_G;R).$$
In particular if $X$ is a point,
$$H^*_G(X;R) = H^*(BG;R).$$
The projection
$$\pi : X_G = X \times _G EG \longrightarrow EG/G = BG$$
defines an action of $H^*(BG;R)$ on $H^*_G(X;R) = H^*(X_G;R)$ making $H^*(X_G;R)$ into a graded module over the graded algebra $H^*(BG;R)$, via
$$(a,b) \in H^*(X_G;R) \times H^*(BG;R) \mapsto a \cup \pi^*(b) \in H^*(X_G;R).$$
This module structure helps measure the isotropy groups.  We will apply equivariant cohomology in the case where $G = O(3)$ and the coefficient group is ${\mathbb Z}_2$, and in this case,
$$H^*(BO(3); {\mathbb Z}_2) = P[w_1,w_2,w_3]$$
a polynomial algebra on generators $w_1$, $w_2$ and $w_3$ of dimensions one, two and three, which we can identify with the Stiefel-Whitney classes of the universal bundle over $BO(3)$.

If $M$ is a compact Riemannian manifold, we let $M_0$ denote the submanifold of ${\mathcal M} = \hbox{Map}(S^2 ,M)$ consisting of constant maps.  Thus
$$M_0 = \{ f \in L^{2\alpha }_1(S^2 ,M) : E_{\alpha }(f) = 0 \}$$
is a submanifold consisting entirely of critical points, which is diffeomorphic to $M$.  We let
$${\mathcal M}^a_\alpha = \{ f \in {\mathcal M}: E_\alpha (f) \leq a \},$$
and use equivariant Morse theory to describe how the relative cohomology
$$H_{O(3)}^*({\mathcal M}^a_\alpha ,M_0,{\mathbb Z}_2)$$
changes as $a$ increases.  We start the process with:

\vskip .1 in
\noindent
{\bf Lemma~8.1.} {\sl There is an $\varepsilon > 0$ such that $M_0$ is a strong deformation retract of ${\mathcal M}^\varepsilon_\alpha $.}

\vskip .1 in
\noindent
This is proven in \cite{SU1}; see Theorem~2.6.  It is not hard to show that $M_0$ is a nondegenerate critical submanifold.

\vskip .1 in
\noindent
We would like to establish Morse inequalities for $E_\alpha $ when $E_\alpha $ is a Morse-Bott function, that is when all of its critical points of positive energy lie on compact nondegenerate critical submanifolds.  It follows quickly from Condition C that only finitely many such submanifolds lie within ${\mathcal M}^a$ for $a < \infty $.

We establish Theorem A and a modification of Theorem B as in Bott's presentation \cite{Bo}.

\vskip .1 in
\noindent
{\bf Theorem~A.} {\sl If the interval $[a,b]$ contains no critical values for $E_\alpha $, then ${\mathcal M}^a_\alpha $ is a strong deformation retract of ${\mathcal M}^b_\alpha $.}

\vskip .1 in
\noindent
Theorem~A is proven by an argument virtually identical to the finite-dimensional case treated in \cite{Mil1}; the deformation retraction is constructed by following the flow lines of $- {\mathcal X}$, where ${\mathcal X}$ is the gradient-like vector field constructed before.  We say that $c$ is a critical value for $F$ if there exists a critical point $f$ for $F$ with $F(f) = c$; the following theorem is proven in \cite{U}, for a Morse perturbation $F$ of $E_{\alpha }$:

\vskip .1 in
\noindent
{\bf Theorem~B.} {\sl If $[a,b]$ contains a single critical value $c$ and there is exactly one Morse nondegenerate critical point $p$ such that $F(p) = c$, then ${\mathcal M}^b_\alpha $ is homotopy equivalent to ${\mathcal M}^a_\alpha $ with a handle of index $\lambda $ attached.}

\vskip .1 in
\noindent
But we need a version of Theorem~B which applies to $F = E_\alpha $, which has nondegenerate critical submanifolds instead of nondegenerate critical points.  For a nondegenerate critical submanifold $K \subset {\mathcal M}$, the normal bundle to $K$ is defined via the inner product $\langle \langle \cdot , \cdot \rangle \rangle$ defined by (\ref{E:innerproduct}):  If $f \in K$,
$$\nu K_f = \{ V \in T_f{\mathcal M} : \langle\langle V,W \rangle\rangle = 0 \hbox{ for all } W\in T_f K \}.$$
Referring back to the operator $A$ defined by (\ref{E:canonicalfredholm}), we see that the normal bundle has a finite-dimensional subbundle $\nu _-K$ generated by eigenvectors corresponding to the negative eigenvalues of $A$.  The rank of the bundle $\nu _-N$ is the {\em Morse index\/} of the nondegenerate critical submanifold $K$.  Let
$$D(\nu _-K) = \{ V \in \nu _-K : \| V \| \leq 1 \}, \qquad S(\nu _-K) = \{ V \in \nu _-K : \| V \| = 1 \},$$
the unit disk and unit sphere bundles in the negative normal bundle $\nu _-K$.  Then Theorem~B has the following generalization:

\vskip .1 in
\noindent
{\bf Theorem B*.} {\sl If $[a,b]$ contains a single critical value $c$ for $F$, the set of critical points with value $c$ forming a nondegenerate critical submanifold $K$ of Morse index $\lambda $, then ${\mathcal M}^b_\alpha $ is homotopy equivalent to ${\mathcal M}^a_\alpha $ with the disk bundle $D(\nu _-K)$ attached to ${\mathcal M}^a_\alpha $ along $S(\nu _-K)$.}

\vskip .1 in
\noindent
The proof is a straightforward modification of the argument for Theorem~B.

The above theorem gives an isomorphism on cohomology
$$H^k({\mathcal M}^b_\alpha ,{\mathcal M}^a_\alpha ;{\mathbb Z}) \cong H^k(D(\nu _-K),S(\nu _-K);{\mathbb Z}).$$
When the bundle $\nu _-K$ is orientable and the group $G$ preserves orientation of the normal bundle, we can use the Thom isomorphism theorem to evaluate the right-hand side:
$$H^k(D(\nu _-K),S(\nu _-K);{\mathbb Z}) \cong H^{k-\lambda }(K;{\mathbb Z}),$$
where $\lambda $ is the Morse index of $K$.  Using Theorem A and Theorem B*, we can now establish equivariant Morse inequalities for ${\mathcal M}^a$ by induction exactly as in Bott's survey article \cite{Bo}.

\vskip .1 in
\noindent
We are interested in applying this theory to the case where $M = S^n$, the standard round $n$-sphere with the constant curvature metric, normalized so all sectional curvatures are one.  Within $M$ we have the nondegenerate submanifold $M_0$ of constant maps, as well as the submanifold $M_{4\pi}$ of prime totally geodesic two-spheres $f : S^2 \rightarrow S^n$, a submanifold  in which all elements have Morse index $n-2$.  For $\alpha $ close to one, the elements of $M_{4\pi}$ are also critical points for the $\alpha $-energy $E_\alpha $ and $M_{4\pi}$ is also a nondegenerate critical submanifold for $E_\alpha $ of Morse index $n - 2$ (although the $\alpha $-energy at points of $M_{4\pi}$ is larger than the usual energy, when $\alpha > 1$).  The action of $O(3)$ on $M_{4\pi}$ is free and the quotient is diffeomorphic to the Grassmann manifold $G_3({\mathbb R}^{n+1})$ of three-planes in ${\mathbb R}^{n+1}$.  Moreover the negative normal bundle $\nu_-M_{4\pi}$ is trivial, hence orientable for any coefficient ring $R$.

\vskip .1 in
\noindent
{\bf Proposition~8.2.} {\sl When $M = S^n$ is the standard round $n$-sphere, every critical point for $E$ is either an imbedded totally geodesic two-sphere with Morse index $n-2$, or has Morse index $\geq 2n-4$.}

\vskip .1 in
\noindent
Proof of Proposition~8.2:  We must show that all minimal two-spheres in $S^n$ with the standard round metric are either totally geodesic of area $4\pi$, or have Morse index at least $2n-4$.  We first note that Corollary~7.2 shows that nontrivial branched covers of totally geodesic minimal two-spheres have Morse index at least $2n-4$.

The theory of harmonic maps $f : S^2 \rightarrow S^n$ which are not branched covers of totally geodesic minimal two-spheres was studied by Calabi \cite{Cal} and we recall some of the high points of his seminal article.  The composition of a harmonic map $f : S^2 \rightarrow S^n$ with the inclusion $S^n \subseteq {\mathbb R}^{n+1}$ is a vector-valued map which in deference to Calabi's notation, we write as
$$X : S^2 \longrightarrow {\mathbb R}^{n+1} \quad \hbox{with} \quad X \cdot X = 1.$$
If $z$ is a local holomorphic coordinate centered at $p \in S^2$ we let $E^{0}(p)$ be the subspace of ${\mathbb R}^{n+1}$ generated by $X(p)$, and for $k \in {\mathbb N}$, we let $E^{k}(p)$ denote the subspace of ${\mathbb R}^{n+1}$ spanned by the real and imaginary parts of the vectors
$$X, \quad \frac{\partial X}{\partial z}, \quad \cdots , \quad \frac{\partial ^kX}{\partial z^k}.$$
For a generic choice of $p$, the dimension of $E^{k}(p)$ is $2k+1$, and we obtain an increasing filtration
$$E^{0}(p) \subset E^{1}(p) \subset E^{2}(p) \subset \cdots \subset E^{m}(p) = E^{m+1}(p) =  \cdots ,$$
which stabilizes at some $E^{m}(p)$.  We can next define a family of normal spaces
\begin{multline*} {\bf L}(p) =  E^{0}(p) \cap E^{1}(p)^\bot, \ \  {\bf L_1}(p) =  E^{2}(p) \cap E^{1}(p)^\bot, \\  {\bf L_2}(p) = E^{3}(p) \cap E^{2}(p)^\bot, \cdots, {\bf L_{m-1}}(p) = E^{m}(p) \cap E^{m-1}(p)^\bot,\end{multline*}
and as $p$ varies we obtain a family of line bundles ${\bf L}, {\bf L_1}, {\bf L_2}, \ldots, {\bf L_{m-1}}$.  (It is easy to extend these bundles to the isolated points at which ranks of derivatives drop.)  Calabi proves that
$$\frac{\partial ^kX}{\partial z^k} \quad \hbox{projects to a holomorphic section of} \quad {\bf L_k},$$
when we use the Koszul-Malgrange connection \cite{KoMa}, and since this section is easily checked to have a zero of order $2k +2 $ at the point $\infty \in S^2$, we see that
\begin{equation} c_1({\bf L_k})[S^2] \geq 2k + 2, \quad \hbox{for} \quad 1 \leq k \leq m-1. \label{E:cilk}\end{equation}
Calabi shows that the $(2m+1)$-dimensional space generated by
$$X, \quad \frac{\partial X}{\partial z}, \quad \cdots , \quad \frac{\partial ^mX}{\partial z^m}, \quad \frac{\partial X}{\partial \bar z}, \quad \cdots , \quad \frac{\partial ^mX}{\partial \bar z^m}$$
is preserved by parallel translation along the two-sphere.  If we denote this subspace of ${\mathbb R}^{n+1}$ by ${\mathbb R}^{2m+1}$, we have a factorization
$$X : S^2 \longrightarrow S^n \cap {\mathbb R}^{2m+1} \subseteq S^n \subseteq {\mathbb R}^{n+1},$$
and we say that $X$ is {\em full\/} in ${\mathbb R}^{2m+1}$.  Calabi shows that the area of $f(S^2)$ is $2m(m+1)\pi \geq 12\pi$, thus it is at least three times the area of the totally geodesic two-spheres.

Following the proof of Theorem~5.1, we can extend our sequence of holomorphic line bundles to achieve a holomorphic direct sum decomposition of the normal bundle ${\bf N}$ to $f$:
$${\bf N} = {\bf L}_1 \oplus {\bf L}_2 \oplus \cdots \oplus {\bf L}_{n-2}, \ \  \hbox{where} \ \ c_1({\bf L}_i)[S^2] = - c_1({\bf L}_{n-i-1})[S^2],$$
with the Chern classes of the first $m-1$ summand being positive, as given by (\ref{E:cilk}), and the last $m-1$ summands being negative.  This leaves $n-2m$ remaining summands which must be flat with the trivial parallel connection.

We can now obtain an estimate on the Morse index of $f$ by adding together contributions from the various summands.  We claim that each trivial summand of ${\bf N}$ contributes two to the Morse index of $f$.  For a trivial summand of the normal bundle ${\bf N}$ it follows from Lemma~5.2 that we can replace the metric by one of constant length, and let $E_3$ denote the canonical unit-length section.  The second variation for energy is the same as that for area $A$ for minimal two-spheres in $S^n$, so we can apply the formula of Proposition~5.3.4 of \cite{Mo2}, and this reduces to the formula
\begin{equation} d^2A(f) (gE_3,hE_3) = \int _{S^2} \langle (\Delta g - 2g) h dA \quad \hbox{where} \quad g,h : S^2 \rightarrow {\mathbb R} \label{E:indexareatwosphere}\end{equation}
are smooth functions and $\Delta $ is the standard Laplace operator for the induced metric on $S^2$.  It follows from the Hersch estimate (Theorem~7.6 of \cite{CM}) that the first nonzero eigenvalue of the Laplace operator for the metric on $S^2$ induced by $f : S^2 \rightarrow S^n$ is $\leq 8\pi/A$, where $A$ is the area of the induced metric.  Since this area is at least $12 \pi$, the second variation formula (\ref{E:indexareatwosphere}) shows that each flat summand of ${\bf N}$ contributes at least two to the Morse index of $f$.

On the other hand, holomorphic sections of the first $m-1$ summands or antiholomorphic sections of the last $m-1$ summands contribute to the Morse index of $f$, and the dimension of the space of holomorphic sections of a holomorphic line bundle over $S^2$ is one greater than the first Chern class.  Thus we get the estimate
\begin{multline*} \hbox{Morse index of $f$} \geq 2n - 4m +  2(5 + 7 + \cdots + 2m+ 1) \\ = 2n - 4m +  2(m+1)^2 - 8 \geq 2n-4,\end{multline*}
completing the proof of the Lemma.   The above argument reproves the estimates of Ejiri \cite{E} for harmonic two-spheres in $S^n$ which are full within ${\mathbb R}^{n+1}$.  QED

\vskip .1 in
\noindent
Recall that when $M = S^n$ is given the round metric of constant curvature one, we let $M_0$ denotes the submanifold of constant maps, $M_{4\pi}$ the submanifold of totally geodesic imbeddings $f : S^2 \rightarrow M$.

\vskip .1 in
\noindent
{\bf Theorem~8.3.} {\sl We have isomorphisms
\begin{multline*} H_{O(3)}^k({\mathcal M} ,M_0,R) \cong H_{O(3)}^{k-n+2}(M_{4\pi};R) \\ \cong H^{k-n+2}(G_3({\mathbb R}^{n+1});R), \quad \hbox{for $k \leq 2n-6$.}\end{multline*}}

\noindent
Proof:  The proof of the previous Proposition shows that all critical points of $E : {\mathcal M} - {\mathcal M} ^{8\pi-} \rightarrow {\mathbb R}$ have Morse index at least $2n-4$.  It is also true that for $\alpha  > 1$ sufficiently small, there is no critical point of  $E_\alpha  : {\mathcal M} - {\mathcal M} ^{8\pi-} \rightarrow {\mathbb R}$ of Morse index $< 2n-3$, because it either converges to a harmonic map which has index at least $2n-4$ or bubbles into at least two two-spheres each having Morse index at least $n-2$.  Hence it follows from equivariant Morse theory that when $\varepsilon > 0$ is sufficiently small,
$$H_{O(3)}^k({\mathcal M} ,{\mathcal M} ^{8\pi-\varepsilon }_\alpha ,R) = 0,  \quad \hbox{for $k \leq 2n-5$.}$$
It therefore follows from the long exact sequence of the triple $({\mathcal M} ,{\mathcal M}^{8\pi-\varepsilon }_\alpha ,M_0)$ that the map induced by inclusion
\begin{equation} H_{O(3)}^k({\mathcal M} ,M_0,R) \longrightarrow H_{O(3)}^k({\mathcal M}^{8\pi-\varepsilon }_\alpha ,M_0,R) \label{E:injective}\end{equation}
is an isomorphism when $k \leq 2n-6$ and injective when $k = 2n-5$.  This implies Theorem~8.3, since by Theorem B*,
$$H_{O(3)}^k({\mathcal M}^{8\pi-\varepsilon }_\alpha ,M_0,R) \cong  H_{O(3)}^{k-n+2}(M_{4\pi};R) \\ \cong H^{k-n+2}(G_3({\mathbb R}^{n+1});R),$$
for all $k \in {\mathbb N} \cup \{0\}$.

\vskip .1 in
\noindent
This calculation shows that the action of $H^*(BO(3); R)$ on $H_{O(3)}^k({\mathcal M},M_0,R)$ is trivial for $k \leq 2n-6$.  Although we have calculated the low-dimensional equivariant cohomology using the constant curvature metric, the cohomology is of course independent of the choice of Riemannian metric on $M = S^n$.  

In particular, when $S^n$ has the metric of constant curvature one, the results of Calabi show that there are no minimal two-spheres $f : S^2 \rightarrow S^n$ of Morse index $\leq 2n-5$ which cover minimal projective planes $f_0 :{\mathbb R}P^2 \rightarrow S^n$.  The minimal projective planes of lowest area and lowest Morse index in $S^n$ are the constant curvature Veronese surfaces $f : S^2 \rightarrow S^4 \subseteq S^n$ as described in Chern \cite{Ch}, which have Morse index $10$ in $S^4$ and also have Morse index $\geq 2n-4$ (by the proof of Lemma~8.2) when considered as lying within $S^n$ for $n \geq 5$.

\section{The Morse-Witten complex}
\label{S:morsewitten}

Now that we know the equivariant cohomology of low degree, we want to use this to prove existence of minimal two-spheres on Morse index $\leq 2n-6$ in manifolds which satisfy our curvature hypothesis (\ref{E:1/3}).  The modern point of view is to organize this via a Morse-Witten complex, and we now review that idea.

Suppose first that $M$ is a finite-dimensional smooth manifold, that $F : M \rightarrow [0, \infty)$ is a proper Morse function and that $X$ is a \lq\lq gradient-like" vector field for $F$, as defined in \S 3 of \cite{Mil2}.  (For example, we might take $X$ to be the gradient of $F$ with respect to a Riemannian metric on $M$.)  We choose an orientation for each critical point, that is, we choose an orientation of its unstable manifold, and let $C_m(F,-X)$ denote the free abelian group generated by the oriented critical points $\hbox{Crit}_m(F)$ of $F$ of Morse index $m$.  If we choose $X$ so that the stable and unstable manifolds of its critical points have transverse intersection (which is true for a residual set of $X$), we can define a boundary operator
\begin{multline*} \partial : C_m(F,-X) \rightarrow C_{m -1}(F,-X) \quad \hbox{by} \\ \partial (p) = \sum _{q \in \hbox{Crit}_{m-1}(F)} n(p,q) q, \quad \hbox{for $p \in \hbox{Crit}_m(F)$,} \end{multline*}
where $n(p,q)$ is the signed number of trajectories for $-X$ from $x$ to $y$.  Construction of this boundary is implicit in Milnor's proof of the h-cobordism theorem \cite{Mil2}, and the theorems of standard Morse theory show that $\partial \circ \partial = 0$, and that the resulting chain complex
\begin{equation} \longrightarrow C_{m+1}(F,-X) \longrightarrow C_m(F,-X) \longrightarrow C_{m -1}(F,-X) \longrightarrow \label{E:MorseWitten1}\end{equation}
computes the homology of $M$.  Of course, we can dualize and construct a corresponding cochain complex
\begin{multline} \longrightarrow C^{m-1}(F,-X) \longrightarrow C^m(F,-X) \longrightarrow C^{m +1}(F,-X) \longrightarrow , \\ \hbox{where} \quad C^m(F,-X) = \hbox{Hom}(C_m(F,-X),{\mathbb Z}). \label{E:MorseWitten2}\end{multline}
We will call either $C_*(F,-X)$ or $C^*(F,-X)$ the {\em Morse-Witten complex} for $F$ and $X$.  More details of the resulting Morse homology can be found in the excellent treatment by Schwartz \cite{Sch}, although in the finite-dimensional case the analysis can be replaced by easier arguments of Milnor \cite{Mil2}.

We can use the Morse theory on Banach manifolds of \S~\ref{S:MorseBanach} to extend the construction of the Morse-Witten complex in \cite{Mil2} to a perturbation of $E_{\alpha}$ which has nondegenerate critical points.  According to Theorem~4.5.3 of \cite{Mo2}, we can take the perturbation to be
\begin{equation} E_{\alpha ,\psi} (f) = E_{\alpha } (f) + \int _{S^2} f \cdot \psi dA, \label{E:Morseperturbation}\end{equation}
where $\psi$ is a suitable element of $L^2_k(S^2,{\mathbb R}^N)$.  Condition C still holds for this perturbation within the completion $L^{2\alpha }_1(S^2,M)$, and implies that there are only finitely many critical points for $E_{\alpha ,\psi}$ within
$${\mathcal M}^{a-}_{\alpha ,\psi} = \{ f \in L^{2}_k(S^2,M) : E_{\alpha ,\psi}(f) < a \},$$
for any choice of $a \in {\mathbb R}$.  We can then choose the gradient-like vector field ${\mathcal X}$ so that stable and unstable manifolds have transverse intersection, which enables us to define the Morse-Witten complex as in the finite-dimensional case.

But we can also consider the case of a perturbation of $E_\alpha $ which is equivariant Morse, all nonconstant critical points lying on isolated nondegenerate orbits for the symmetry group, in our case either $\hat G = O(3)$ or its identity component $G = SO(3)$.  By Theorem~2.6 in \cite{SU1} there is a small $\varepsilon > 0$ such that
$${\mathcal M}^{\varepsilon}_{\alpha ,\psi} = L^{2}_k(S^2,M)^\varepsilon = \{ f \in L^{2}_k(S^2,M) : E_\alpha (f) \leq \varepsilon \}$$
possesses the space $M_0$ of constant maps as a strong deformation retract.

\vskip .1 in
\noindent
{\bf Theorem~9.1.} {\sl The equivariant function $E_\alpha $, restricted to ${\mathcal M}^{a-}_\alpha - {\mathcal M}^\varepsilon_\alpha $ for any $a > 0$ and a suitably small $\varepsilon > 0$, can be perturbed to an $SO(3)$-equivariant Morse function which is $C^2$ close to $E_\alpha $ and satisfies Condition C on the completion $L^{2\alpha }_1(S^2,M)$.}

\vskip .1 in
\noindent
Proof:  This would follow from the Density Lemma~4.8 from Wasserman's treatment of equivariant Morse theory \cite{Wa}, except that our $G$-manifolds are not compact, and we want to ensure that our perturbed functions have smooth critical points and satisfy Condition C on the completion $L^{2\alpha }_1(S^2,M)$.  We will explain how to modify Wasserman's argument so that it applies to our situation.  We will assume known the theory of slices and tubes from Chapter~2 of \cite{DK}, which carries over to infinite-dimensional $G$-manifolds with almost no change.

Note that by Lemma~3.1 we can assume that the isotropy group of $E_{\alpha }$ at each critical point is finite, and the isotropy is also finite on a small neighborhood of the compact critical locus for $E_\alpha $ within the open set ${\mathcal M}^{a-} - {\mathcal M}^\varepsilon$.  If we make a $C^2$ small perturbation, the critical locus will remain in this neighborhood.

We give ${\mathcal M}^{a-}_\alpha - {\mathcal M}^\varepsilon_\alpha$ an increasing filtration by open sets
$${\mathcal U}_n = \{ f \in {\mathcal M}^{a-}_\alpha - {\mathcal M}^\varepsilon_\alpha : \hbox{ the isotropy group $G_{f}$ has order $\leq n$ } \},$$
where $n \in {\mathbb N} \cup \{ \infty \}$.  We will prove by induction on $n$ that for $n \in {\mathbb N}$ the restriction of $E_\alpha $ to ${\mathcal U}_n$ can be perturbed to a $G$-equivariant Morse function.  We will do this by constructing a family of perturbed functions,
$$E_{\alpha ,F} (f) = E_{\alpha} (f) + F(f), \quad \hbox{for} \quad F \in {\mathcal F}_n \subseteq C^2(L^{2}_k(S^2,M),{\mathbb R}),$$
to which we can apply the Baire category theorem.  In fact the perturbations ${\mathcal F}_n$ will be a uniformly closed subspace of the space of $C^2$ functions
$$F: L^2_k(S^2,M) \longrightarrow {\mathbb R}$$
with the {\em strong topology\/} as described in \S4 of Chapter 2 of \cite{Hir}.  As a basis for the open sets of this topology, we can take
$${\mathcal N}(F,\varepsilon) = \{ G \in C^2(L^2_k(S^2,M),{\mathbb R}) : d((j_2F)(f),(j_2G)(f))  < \varepsilon (f) \},$$
where $j_2$ denotes the two-jet and $\varepsilon : L^2_k(S^2,M) \rightarrow (0,1)$ is an arbitrary positive continuous function.  This is a Baire space by Theorem~4.2 of Chapter 2 of \cite{Hir}, so that the intersection of countably many open dense subsets of ${\mathcal F}_n$ is dense within ${\mathcal F}_n$.

\vskip .1in
\noindent
{\bf Start of induction $n = 1$:}  The action of $G$ on ${\mathcal U}_1$ is free and just as in the finite-dimensional case the quotient ${\mathcal U}_1/G$ is a smooth manifold, in fact, ${\mathcal U}_1$ is a principal $G$-bundle over ${\mathcal U}_1/G$.  Since $L^{2}_k(S^2,M)$ has $C^2$ partitions of unity, we can construct countable open covers $\{ U_{1,i} : i \in {\mathbb N} \}$ and $\{ V_{1,i} : i \in {\mathbb N} \}$ of ${\mathcal U}_1$ consisting of $G$-invariant subsets with $\bar V_{1,i} \subseteq U_{1,i}$, together with a $G$-invariant partition of unity $\{ \eta _{1,i} : i \in {\mathbb N} \}$ subordinate to $\{ U_{1,i} : i \in {\mathbb N} \}$ with the support of $\eta _{1,i}$ contained within $\bar V_{1,i}$.  Moreover, we can arrange that each $U_{1,i}$ contains a local slice ${\mathcal S}_{f_{1,i}}$ centered at $f_{1,i} \in U_{1,i}$ such that each $G$-orbit in $U_{1,i}$ intersects ${\mathcal S}_{f_{1,i}}$ exactly once.  We can then define a family of perturbed functions
\begin{multline} E_{\alpha ,F} (f) = E_{\alpha } (f) + \sum _i \eta _{1,i} F_{1,i}(f), \\ \hbox{where} \quad F_{1,i}(f) = \int _{S^2} f \cdot \psi _{1,i} dA, \quad \hbox{for $f \in {\mathcal S}_{f_{1,i}}$,} \label{E:Morseperturbation1}\end{multline}
with $\psi _{1,i} \in L^2_k(S^2,{\mathbb R}^N)$ and $F_{1,i}$ is extended to be $G$-invariant.  We let
$${\mathcal F}_1 = \left\{ \sum _i \eta _{1,i} F_{1,i} \right\},$$
a family of functions on $L^2_k(S^2,M)$, where the open covers $\{ U_{1,i} \}$ and $\{ V_{1,i} \}$ are fixed but the partition of unity varies subject to the conditions
\begin{equation}  \hbox{supp}(\eta _{1,i}) \subseteq \bar V_{1,i}, \quad \|  \eta _{1,i} \|_{C^2} \leq C_0, \label{E:conditions1}\end{equation}
where $C_0$ is a given bound.  Then ${\mathcal F}_1$ is closed in the strong topology, and the usual argument shows that the restriction of $E_{\alpha ,F}$ to ${\mathcal U}_1$ is an equivariant Morse function for a generic choice of perturbation $F \in {\mathcal F}_1$.

\vskip .1in
\noindent
{\bf The case of general $n > 1$:}  We suppose that ${\mathcal F}_{n-1}$ has been constructed, and we let
$${\mathcal N}_n = \{ f \in {\mathcal U}_n : \hbox{ the isotropy group $G_{f}$ has order $n$ } \},$$
a smooth submanifold with ${\mathcal U}_n = {\mathcal U}_{n-1} \cup {\mathcal N}_n$.

We construct countable open covers $\{ U_{n,i} : i \in {\mathbb N} \}$ and $\{ V_{n,i} : i \in {\mathbb N} \}$ of ${\mathcal N}_n$, consisting of $G$-invariant open subsets of ${\mathcal U}_n$ with $\bar V_{n,i} \subseteq U_{n,i}$, together with a $G$-invariant partition of unity $\{ \eta _{n,i} : i \in {\mathbb N} \}$ subordinate to $\{ U_{n,i} : i \in {\mathbb N} \}$ with $\hbox{supp}(\eta _{n,i})$ contained within $\bar V_{n,i}$.  We assume that $U_{n,i} $ contains $f_{n,i} \in {\mathcal N}_n$, and that $U_{n,i}$ possesses a local slice ${\mathcal S}_{f_{n,i}}$ centered at $f_{n,i}$, with each $G$-orbit within $U_{n,i}$ intersecting ${\mathcal S}_{f_{n,i}}$.  We can suppose that ${\mathcal S}_{f_{n,i}}$ is the domain of a tangent-space-valued chart centered at $f_{n,i}$,
\begin{multline*} \phi _{n,i} : {\mathcal S}_{f_{n,i}} \longrightarrow E \oplus F = T_{f_{n,i}}{\mathcal S}_{f_{n,i}}, \\ \hbox{where} \quad \phi _{n,i} ({\mathcal N}_n \cap {\mathcal S}_{f_{n,i}}) = ( E \oplus \{ 0 \} ) \cap \phi _{n,i} ({\mathcal S}_{f_{n,i}}).\end{multline*}
Indeed, according to Bochner's Theorem (see 2.2.1 in \cite{DK}) we can assume that the differential of $\phi _{n,i}$ is the identity, and that $\phi _{n,i}$ takes the action of the isotropy group $G_{f_{n,i}}$ on $f_{n,i}$ to an orthogonal linear action on the tangent space which preserves the direct sum decomposition $E \oplus F$ and is trivial on $E$.  We can assume also that $\phi _{n,i} ({\mathcal S}_{f_{n,i}})$ is a product neighborhood in $E \oplus F$.  

We can think of $\phi _{n,i}$ as defining an $E$-valued chart on ${\mathcal N}_n \cap {\mathcal S}_{f_{n,i}}$, and as in the case $n = 1$ we can define a family of perturbed functions on ${\mathcal N}_n$,
\begin{multline} E_{\alpha ,F} (f) = E_{\alpha } (f) + \sum _i \eta _{n,i} F_{n,i}(f), \\ \hbox{where} \quad F_{n,i}(f) =  \int _{S^2} \phi _{n,i}(f) \cdot \psi _{n,i} dA, \quad \hbox{for $f \in {\mathcal S}_{f_{n,i}}$,} \label{E:Morseperturbation2}\end{multline}
with $\psi _{n,i} \in E$, where $E$ is a linear subspace of $L^2_k(S^2,{\mathbb R}^N)$ and $F_{n,i}$ is extended to be $G$-invariant on $U_{n,i}$.  As in the case $n = 1$, the restriction to $U_{n,i} \cap {\mathcal N}_n$ will be Morse for a generic choice of $\psi _{n,i}$'s.

Finally, we let $P_E$ and $P_F$ denote orthogonal projections to $E$ and $F$ respectively, and define
$$r_i^2 : U_{n,i} \rightarrow [0, \infty) \quad \hbox{by} \quad r_i^2(f) = \| P_F(\phi _{n,i}(f)) \|^2.$$
We then modify (\ref{E:Morseperturbation2}) to the family of perturbations
\begin{multline} E_{\alpha ,F} (f) = E_{\alpha } (f) + \sum _i \eta _{n,i} \left[ \varepsilon _i r_i^2(f) + F_{n,i}(f) \right], \\ \hbox{where} \quad F_{n,i}(f) = \int _{S^2} P_E(\phi _{n,i}(f)) \cdot \psi _{n,i} dA, \quad \hbox{for $f \in {\mathcal S}_{f_{n,i}}$,} \label{E:Morseperturbation3}\end{multline}
Since $r_i^2$ is invariant under the action of the isotropy group $G_{f_{n,i}}$, we can extend $E_{\alpha ,F}$ to $U_{n,i}$ so that it is $G$-invariant.  For a generic choice of $\psi _{n,i}$'s, this will have the property that all critical points within $U_{n,i} \cap {\mathcal N}_n$ are nondegenerate.

To complete the inductive step we now set
$${\mathcal F}_n = {\mathcal F}_{n-1} \cup \left\{ \sum _i \eta _{n,i} \left[ \varepsilon _i r_i^2(f) + F_{n,i}(f) \right] \right\},$$
where the open covers $\{ U_{n,i} \}$ and $\{ V_{n,i} \}$ are fixed but the partition of unity varies subject to the conditions
\begin{equation}  \hbox{supp}(\eta _{n,i}) \subseteq \bar V_{n,i}, \quad \|  \eta _{n,i} \|_{C^2} \leq C_0, \label{E:conditions2}\end{equation}
where $C_0$ is a given bound.  Once again, the Baire category argument implies that the resulting perturbed functions will be Morse for generic elements of ${\mathcal F}_n$.

Finally, we observe that our perturbations do not involve the derivatives of the map $f$ and can be made arbitrarily small.  The fact that (\ref{E:conditions2}) holds for all $n \in {\mathbb N}$ makes it easy to check that Condition C still holds on the completion, and that the regularity argument of Sacks and Uhlenbeck implies that any critical point must be in $L^2_k(S^2,M)$.  These results are discussed further in \cite{Mo2}.  This finishes the proof.

\vskip .1 in
\noindent
Note that if $M$ is $C^\infty$, we can take $k$ to be arbitrary in Theorem~9.1, and thus we can arrange that all critical points are $C^\infty $.

We want to apply Theorem~9.1 to construct an equivariant Morse-Witten complex.  When the metric on $M$ is generic, all nonconstant minimal two-spheres of Morse index $\leq 2n-5$ in $M$ are immersions and lie on Morse nondegenerate critical submanifolds, each a $\widehat {\mathcal G}$-orbit, where $\widehat {\mathcal G}$ is the group of conformal and anticonformal automorphisms of $S^2$.  There are only finitely many of these by Theorem~6.2 and according to Theorem~2.2, these perturb to finitely many $G$-orbits of $\alpha $-energy critical points, each a nondegenerate critical submanifold.  Let $U$ be an open subset which contains all critical points of Morse index $\leq 2n-5$, and suppose that $V$ is a larger open set such that $U \subseteq \bar U \subseteq V$ and $V$ contains no additional critical points.  Since $L^{2\alpha }_1(S^2,M)$ has $C^2$ partitions of unity, we can construct a smooth function $\eta : L^{2\alpha }_1(S^2,M) \rightarrow [0,1]$ which is zero on $U$ and one on $L^{2\alpha }_1(S^2,M) - V$.  Then we let
$$F : L^{2\alpha }_1(S^2,M) \longrightarrow {\mathbb R} \quad \hbox{by} \quad F = (1- \eta )E_\alpha + \eta \tilde F,$$
where $\tilde F$ is an equivariant Morse approximation to $E_\alpha $ as described in the previous paragraph, all of its critical points lying on Morse nondegenerate critical submanifolds, each an orbit for $G$.  This $F$ is an equivariant Morse function which agrees with $E_\alpha $ on $U$.

Note that except for points lying within the space $M_0$ of constant maps, the isotropy groups of $E_\alpha $ are always finite by Lemma~3.1, and we can arrange that $F$ satisfies the same condition.  Let $\{ \phi _t : t \in [0,\infty) \}$ be the one-parameter semigroup of diffeomorphisms of ${\mathcal M}$ corresponding to $- {\mathcal X}$, this semigroup being defined for all positive $t$ because $- {\mathcal X}$ is positively complete.  Note that the isotropy group for the action of either $\widehat G$ or $G$ is constant along any orbit for $- {\mathcal X}$.  Recall that the {\em unstable manifold\/} $W_f(F,{\mathcal X})$ of a critical point $f \in {\mathcal M}$ consists of the images of all trajectories $t \mapsto \phi_t(g)$ which start at $f$ in the remote past, in other words, such that $\phi_t(g) \rightarrow f$ as $t \rightarrow - \infty$.  Similarly, the {\em stable manifold\/} $W_f^*(F,{\mathcal X})$ of a critical point $p$ consists of the images of trajectories $t \mapsto \phi_t(g)$ such that $\phi_t(g) \rightarrow f$ as $t \rightarrow \infty$.    

\vskip .1 in
\noindent
{\bf Lemma~9.2.} {\sl The unstable and stable manifolds $W_f(F,{\mathcal X})$ and $W_f^*(F,{\mathcal X})$ are indeed submanifolds of ${\mathcal M}$.}

\vskip .1 in
\noindent
The proof follows from the explicit form of the gradient-like vector field ${\mathcal X}$ near the critical point $f$ and the properties of local flows for vector fields on Banach manifolds.

\vskip .1 in
\noindent
It follows from Condition C and the explicit form of ${\mathcal X}$ that any orbit $g \mapsto \phi _t(g)$ converges to some critical point as $t \rightarrow \infty $ (although it may come close to several other critical points first).

\vskip .1 in
\noindent
{\bf Lemma~9.3.} {\sl We can adjust the $G$-equivariant gradient-like vector field ${\mathcal X}$ so that if $\lambda _f$ is the Morse index of the critical point $f \in {\mathcal M}$ and $\lambda _g$ is the Morse index of the critical point $g$,
\begin{enumerate}
\item $\lambda _f \leq \lambda _g \Rightarrow W_f({\mathcal X})\cap W_g^*({\mathcal X})$ is empty, while
\item $\lambda _f > \lambda _g \Rightarrow W_f({\mathcal X})\cap W_g^*({\mathcal X})$ is a submanifold of dimension $\lambda _f - \lambda _g + 3$.
\end{enumerate}}

\vskip .1 in
\noindent
Proof:  If $W_p({\mathcal X}) \cap W_g^*({\mathcal X})$ is nonempty, then $F(f) > F(g)$ and we can choose a regular value $c$ for $F$ such that $F(f) > c > F(g)$, noting that we can choose $c$ as close to $F(f)$ as we want.  Then ${\mathcal N} = F^{-1}(c)$ is a codimension one submanifold of ${\mathcal M}$, and we let ${\mathcal S}_h$ be a section for the $G$-action through the point $h \in {\mathcal N} $, a submanifold of codimension three.  Note that
$$\hbox{dim } ({\mathcal N} \cap {\mathcal S}_h \cap W_f({\mathcal X})) = \lambda _f - 1, \qquad \hbox{codim } ({\mathcal N} \cap W_g^*({\mathcal X})) = \lambda _g.$$

If we choose $c$ sufficiently close to the Morse nondegenerate critical point $p$, we see that $S = {\mathcal N} \cap {\mathcal S}_h \cap W_f({\mathcal X})$ is a $(\lambda _f - 1)$-dimensional sphere and lies in an open tubular neighborhood $U$ of ${\mathcal N}$ with diffeomorphism
$\psi : U \rightarrow S \times V$, where $V$ is an open ball in a Banach space, with $\psi (S) = S \times \{ 0 \}$.  Let
$$N = U \cap  W_h^*({\mathcal X}), \quad \hbox{and consider} \quad H = \pi \circ \psi : N \rightarrow V,$$
where $\pi $ is the projection on the second factor.  We can check that $H$ is a Fredholm map with Fredholm index $\lambda _f - 1 - \lambda _g $, and
$$\psi ^{-1}(S \times \{x\})\cap N = \{ s \in N : H(s) = x \}.$$
Assuming that ${\mathcal X}$ (and hence $N$) is sufficiently smooth, we can use the Sard-Smale Theorem of \cite{Sm2} to choose a regular value $x$ for $H$.  Finally, we can choose $x$ as close as we want to $0$, and construct an isotopy from a neighborhood of $S$ to itself which carries $S \times \{0\}$ to $S \times \{x\}$.  Since the conditions defining pseudogradient are open, we can replace ${\mathcal X}$ by a new gradient-like vector field so that
$${\mathcal N} \cap {\mathcal S}_f \cap W_f({\mathcal X}) \quad \hbox{and} \quad {\mathcal N} \cap W_g^*({\mathcal X})$$
have transverse intersection.  Then the intersection of
$${\mathcal N} \cap W_f({\mathcal X}) \quad \hbox{and} \quad {\mathcal N} \cap W_g^*({\mathcal X})$$
will be a submanifold of dimension $\lambda _f - \lambda _g + 3$.

\vskip .1 in
\noindent
Since $F$ satisfies Condition C, there are only finitely many $G$-orbits of critical points within
$${\mathcal M}^a_\alpha = \{ f \in L^{2\alpha }_1(S^2,M) : F (f) \leq a \},$$
for any choice of $a \in {\mathbb R}$.  Lemma 9.2 implies that after perturbation of $- {\mathcal X}$ unstable manifolds of critical orbits can only intersect stable manifolds of critical points of lower index.  Using the arguments in \cite{Mil2}, we can then modify the function $F$ to a new function $\tilde F$, without changing $G$-orbits of critical points or trajectories between critical points, so that $\tilde F$ is a self-indexing equivariant Morse function, one for which each critical point $f \in {\mathcal M}^a$ satisfies
$$\tilde F(f) = \hbox{Morse index of $f$}.$$

Now we can use further arguments within \cite{Mil2} to define the Morse-Witten complexes (\ref{E:MorseWitten1}) and (\ref{E:MorseWitten2}).  Let us focus on the cohomology complex.  For $\lambda \in \{ 0 \} \cup {\mathbb N}$, we let
$$C^\lambda (F,-{\mathcal X}) = H^\lambda_{G} ({\mathcal M}^{\lambda + (1/2)}, {\mathcal M}^{\lambda - (1/2)}; {\mathbb Z}),$$
and define
$$d : C^\lambda (F,-{\mathcal X}) \longrightarrow C^{\lambda + 1} (F,-{\mathcal X})$$
to be the coboundary operator in the triple
$$({\mathcal M}^{\lambda + (3/2)}, {\mathcal M}^{\lambda + (1/2)}, {\mathcal M}^{\lambda - (1/2)}).$$
using the fact that the cohomology of the pair $({\mathcal M}^{\lambda + (1/2)}, {\mathcal M}^{\lambda - (1/2)})$ is concentrated in degree $\lambda $, we then obtain a cochain complex
$$\cdots \longrightarrow C^{\lambda-1}(\tilde F,-{\mathcal X}) \longrightarrow C^\lambda(\tilde F,-{\mathcal X}) \longrightarrow C^{\lambda +1}(\tilde F,-{\mathcal X}) \longrightarrow \cdots ,$$
which we call the Morse-Witten complex, in which each of the coboundary maps $d$ is an $H^*(BG)$-module homomorphism.  In the case where the action of $G$ is free, this does indeed count trajectories connecting critical points in the orbit space ${\mathcal M}/G$, and in keeping with the philosophy of equivariant cohomology, we regard this as counting trajectories in ${\mathcal M}/G$ even when the action is only locally free.

\section{Existence of low index minimal two-spheres}
\label{S:existence}

Finally, we collect the preceding ingredients to prove existence of minimal two-spheres of Morse index $\leq 2n-6$ in manifolds which satisfy the curvature conditions (\ref{E:Kr}) or (\ref{E:1/3bis}), thereby proving Theorem~1.3.

Suppose that $M$ is homeomorphic to $S^n$ and that $M$ is given a Riemannian metric of positive half-isotropic sectional curvatures.  Then it follows from Theorem~6.2 that there are only finitely many $\widehat {\mathcal G}$-orbits of prime minimal two spheres $f : S^2 \rightarrow M$ with Morse index $\leq 2n-5$, where $\widehat {\mathcal G} = PSL(2,{\mathbb C}) \cup R \cdot PSL(2,{\mathbb C})$.  Some of these may be invariant under an antipodal map $A$, and hence double cover minimal projective planes.  If $\alpha _0 > 1$ is sufficiently close to one, each of these $\widehat {\mathcal G}$-orbits perturbs to a single nondegenerate $O(3)$-orbit of critical points for $E_\alpha $ for $\alpha \in (1, \alpha _0)$.  Moreover, for Morse index $\leq 2n-5$, there are no sequences of $\alpha _m$-energy critical points with $\alpha _m \rightarrow 1$ which approach bubble trees with more than one bubble by Theorem~1.4.  Indeed, if we have a sequence of such $\alpha _m$-energy critical points with $\alpha _m \rightarrow 1$, a subsequence must either converge to an orbit of prime minimal two-spheres or a branched cover of such a prime minimal two sphere.  Under the curvature condition (\ref{E:1/3}), any such branched cover must be of order two, since otherwise it would have index $\geq 2n-4$ by Theorem~7.1.

After imposing the center of mass condition described in \S \ref{S:perturbations}, parametrized minimal two-spheres $f: S^2 \rightarrow M$ of Morse index $\leq 2n-5$ are partitioned into several types depending upon their isotropy groups for the action of $\widehat G = O(3)$:
\begin{enumerate}
\item We have prime minimal two-spheres for which the isotropy group of the $\widehat G$-action is trivial.  In this case, the $\widehat G$-orbit consists of two components, one with each of two orientations.
\item We may also have prime minimal projective planes for which the isotropy group $H$ of the $\widehat G$-action is isomorphic to ${\mathbb Z}_2$ and is generated by an orientation-reversing element $A \in \widehat G$.  In this case, the $\widehat G$-orbit consists of a single component.  
\item Next we may have nontrival branched covers of prime minimal two-spheres, for which the isotropy group $H$ of $\widehat G$ is a nontrivial discrete subgroup of $G = SO(3)$.  Under our curvature hypothesis, the branching order must be two. 
\item Finally, we might have nontrivial branched covers of prime minimal projective planes.
\end{enumerate}
Since bubbling is ruled out, sequences of $\alpha _m$-energy critical points with $\alpha _m \rightarrow 1$ must converge to a minimal two-sphere which lies in one of these cases.  We claim that the cases can be distinguished within equivariant cohomology.

We now construct a Morse-Witten complex for the function $E_\alpha$ on the pair $({\mathcal M},M_0)$ as in \S~\ref{S:morsewitten}, and determine the contributions given by Theorem B* of \S \ref{S:equivariant} to the cohomology of this complex when all critical $\widehat G$-orbits are nondegenerate.  We claim that the type of critical orbit obtained is encoded in the action of
$$H^*(BO(3);{\mathbb Z}_2) = P[w_1,w_2,w_3]$$
on the equivariant cohomology of the corresponding generator in the Morse-Witten complex.  Suppose we want to calculate the contribution to the Morse-Witten complex of an isolated nondegenerate critical orbit $K$ with isotropy group $H$, so $K$ is diffeomorphic to $\widehat G/H$.  The standard approach for computing the equivariant cohomology of $K$ is to use the spectral sequence of the fibration
\begin{equation} \begin{matrix} K & \longrightarrow & K_{\widehat G} \cr && \pi \downarrow \cr && B\widehat G, \end{matrix} \label{E:spectral}\end{equation}
as explained in Chapter~5 of McCleary \cite{McC}, the equivariant cohomology of $K$ being the usual cohomology of the total space $K_{\widehat G}$.

1. In the case of trivial isotropy, $K = O(3)$ and the homotopy quotient is contractible, and this forces the action of $P[w_1,w_2,w_3]$ to be trivial.

2. In the case of prime minimal projective planes, in which the isotropy group is isomorphic to ${\mathbb Z}_2$ and is generated by an orientation-reversing element $A \in O(3)$, the orbit space is $K = SO(3)$ which has ${\mathbb Z}_2$ cohomology
$$H^*(SO(3);{\mathbb Z}_2) \cong \Lambda( \xi_1,\xi_2), \quad \hbox{where} \quad \deg \xi _i = i,$$
and the spectral sequence shows that $w_1$ acts freely on the homotopy quotient $K_{\widehat G}$, while $w_2$ and $w_3$ act trivially.

3 and 4. For the remaining cases, we consider the isotropy group under the action of $G = SO(3)$.  For double covers, the isotropy group is ${\mathbb Z}_2$ and is generated by a rotation, while the orbit space $K = SO(3)/{\mathbb Z}_2$ has fundamental group $\pi _1(K) \cong {\mathbb Z}_4$.  It follows from Poincar\'e duality that
$$H^*(K;{\mathbb Z}_2) \cong \Lambda( \eta_1,\eta_2), \quad \hbox{where} \quad \deg \xi _i = i,$$
once again, but the projection $\pi : SO(3) \rightarrow K$ induces a map on ${\mathbb Z}_2$-cohomology which takes $\eta _1$ to zero.  This implies that in the spectral sequence of (\ref{E:spectral}), $w_2$ is not in the image of transgression, which implies that $w_2$ must act freely on the homotopy quotient in this case.  Since the equivariant cohomology of $K$ is just $H_{SO(3)}(K;{\mathbb Z}_2) \cong P[u]$, a polynomial algebra generated by an element $u$ of degree one, it must be that $w_2$ acts by taking $1$ and $u$ to $1 \cdot w_2 = u^2$ and $uw_2 = u^3$ respectively, and in particular the action of $O(3)$ is nontrivial. 

Now we use Theorem B* of \S \ref{S:equivariant} and the Thom isomorphism theorem to show that each critical orbit yields a generator for
$$H^*_{O(3)}({\mathcal M}, M_0,{\mathbb Z}_2) \quad \hbox{as a module over} \quad P[w_1,w_2,w_3]$$
the degree of the generators being shifted upward by the Morse index of the critical orbit.

We now divide the Morse-Witten complex $C^\lambda (F,-{\mathcal X})$ into a direct sum,
$$C^\lambda (F,-{\mathcal X}) = A^\lambda (F,-{\mathcal X}) \oplus B^\lambda (F,-{\mathcal X}),$$
where $H^*(BO(3);{\mathbb Z}_2) = P[w_1,w_2,w_3]$ acts trivially on $A^\lambda $ and nontrivially on $B^\lambda $, the generators of $A^\lambda $ being prime minimal two-spheres with two possible orientations, the generators of $B^\lambda $ being projective planes and double covers of prime minimal two-spheres or projective planes.  Note that
$$\hbox{$P[w_1,w_2,w_3]$ acts trivially on $A^\lambda $} \quad \Rightarrow \quad d(A^\lambda ) \subseteq A^\lambda ,$$
because the differential $d$ commutes with the action of $P[w_1,w_2,w_3]$.  Thus the $A^\lambda $'s form a subcomplex of the Morse-Witten complex, and we have a short exact sequence
$$0 \longrightarrow A^* \longrightarrow C^* \longrightarrow B^* \longrightarrow 0.$$

But we have noted before that in the constant curvature case, there are no double covers of real projective planes with Morse index $\leq 2n-6$, so we have
$$H^k(A^*;{\mathbb Z}_2) \cong H^k(C^*;{\mathbb Z}_2) \cong H^{k-n+2}(G_3({\mathbb R}^{n+1});{\mathbb Z}_2),$$
when $k \leq 2n-6$, the last isomorphism following from Theorem~8.3.  When the metric on $M$ is generic, the Morse inequalites for the subcomplex $A^*$ now show that the number of minimal two-spheres which do not cover projective planes of Morse index $\lambda $ for $n-2 \leq \lambda \leq 2n - 6$ must be at least as large as the number $p_{3}(\lambda - n + 2)$, where $p_{3}(k)$ is the number of $k$-cells in the Schubert cell decomposition for $G_3({\mathbb R}^{n+1})$.  This proves Theorem~1.3 when $\lambda $ lies in the range  $n-2 \leq \lambda \leq 2n - 6$.

Theorem~8.3 does not extend to the case where $ \lambda = 2n - 5$,  because (\ref{E:injective}) is no longer an isomorphism, but we can still extend the proof for Theorem~1.3 to this case:  Since the action of $H^*(BO(3);{\mathbb Z}_2)$ on $H_{O(3)}^{2n-5}({\mathcal M} ,M_0,R)$ is trivial, there is a subspace of $H_{O(3)}^{2n-5}({\mathcal M} ,{\mathcal M} ^{8\pi-\varepsilon }_\alpha ,R)$ of dimension $p_{3}(n - 3)$ on which the action of $H^*(BO(3);{\mathbb Z}_2)$ is also trivial.  We can use this fact, together with the previous argument, to show that if $E_{\alpha , \psi}$ is a perturbation of $E$ with nondegenerate critical points,
$$\hbox{the number of critical points of $E_{\alpha , \psi}$ of Morse index $2n-5$} \geq p_{3}(n - 3).$$
Once again the Nonbubbling Theorem implies that bubbling cannot occur, so when $\alpha \rightarrow 1$ we obtain at least $p_{3}(n - 3)$ minimal two-spheres which cannot cover projective planes for generic Riemannian metrics on $M$.  This finishes the proof of the Main Theorem~1.3 in the remaining case.

\vskip .1 in
\noindent
{\bf Remark~10.1.} The strategy outlined in this section can also be applied to yield a small improvement in the Corollary in \S 1 of \cite{Mo1}:  Instead of considering the full Morse-Witten complex for ${\mathcal M} = \hbox{Map}(S^2 ,M)$ we consider the subcomplex on which $H^*(BO(3);{\mathbb Z}_2) = P[w_1,w_2,w_3]$ acts trivially.  The critical points for this subcomplex cannot double cover minimal projective planes because these critical points must have trivial isotropy group.  We can therefore remove the hypothesis in the Corollary that $M$ have no minimal projective planes of area $< 4 \pi$.

\end{document}